\documentclass[12pt]{article}

\usepackage{graphicx}
\usepackage{latexsym,amsmath,amsfonts,amscd, amsthm}
\usepackage{changebar}
\usepackage{color}
\usepackage{bm}

\newtheoremstyle{plainNoItalics}{}{}{\normalfont}{}{\bfseries}{.}{ }{}

\theoremstyle{plain}
\newtheorem{thm}{Theorem}[section]

\theoremstyle{plainNoItalics}

\newtheorem{defn}[thm]{Definition}
\newtheorem{rem}[thm]{Remark}
\newtheorem{prop}[thm]{Proposition}
\newtheorem{exa}[thm]{Example}

\newcommand{\beq}{\begin{equation}}
\newcommand{\eeq}{\end{equation}}
\newcommand{\beqa}{\begin{eqnarray}}
\newcommand{\eeqa}{\end{eqnarray}}
\newcommand{\bit}{\begin{itemize}}
\newcommand{\eit}{\end{itemize}}
\newcommand{\bedef}{\begin{defn}}
\newcommand{\edefn}{\end{defn}}
\newcommand{\bpro}{\begin{prop}}
\newcommand{\epro}{\end{prop}}

\newcommand{\Dx}{\Delta x}
\newcommand{\Dt}{\Delta t}

\newcommand{\pad}[2]{\frac{\partial #1}{\partial #2}}

\newcommand{\jp}{{j+1/2}}
\newcommand{\jm}{{j-1/2}}

\newcommand{\rf}[1]{(\ref{#1})}

\newcommand{\hF}{\hat{F}}

\def\Box{\mbox{ }\rule[0pt]{1.5ex}{1.5ex}}

\begin{document}

\baselineskip=1.8pc



\begin{center}
{\bf 
A High Order Multi-Dimensional Characteristic Tracing Strategy for 
the Vlasov-Poisson System
}
\end{center}

\vspace{.2in}
\centerline{  
Jing-Mei
Qiu\footnote{Department of Mathematics, University of Houston,
Houston, 77004, USA. E-mail: jingqiu@math.uh.edu. Research supported by
Air Force Office of Scientific Computing YIP grant FA9550-12-0318, NSF grant DMS-1217008 and DMS-1522777.},
Giovanni Russo \footnote{Department of Mathematics and Informatics, University of Catania, Catania, 95125, Italy. Email: russo@dmi.unict.it. 
Research supported by ITN-ETN Marie Curie program 642768.
} 
}

\bigskip
\centerline{\bf Abstract}
\bigskip
In this paper, we consider a finite difference grid-based semi-Lagrangian approach in solving the Vlasov-Poisson (VP) system. 
Many of existing methods are based on dimensional splitting, which decouples the problem into solving linear advection problems, see {\em Cheng and Knorr, Journal of Computational Physics, 22(1976)}.
However, such splitting is subject to the splitting error. If we consider multi-dimensional problems without splitting, difficulty arises in tracing characteristics with high order accuracy. Specifically, the evolution of characteristics is subject to the electric field which is determined globally from the distribution of particle densities via the Poisson's equation. In this paper, we propose a novel strategy of tracing characteristics high order in time via a two-stage multi-derivative prediction-correction approach and by using moment equations of the VP system. With the foot of characteristics being accurately located, we proposed to use weighted essentially non-oscillatory (WENO) interpolation to recover function values between grid points, therefore to update solutions at the next time level. The proposed algorithm does not have time step restriction as Eulerian approach and enjoys high order spatial and temporal accuracy. However, such finite difference algorithm does not enjoy mass conservation; we discuss one possible way of resolving such issue and its potential challenge in numerical stability.
The performance of the proposed schemes are numerically demonstrated via classical test problems such as Landau damping and two stream instabilities. 
\vfill

\newpage

\newpage

\section{Introduction}
\label{sec1}
\setcounter{equation}{0}
\setcounter{figure}{0}
\setcounter{table}{0}

This paper focuses on a high order truly multi-dimensional semi-Lagrangian (SL) approach for the Vlasov-Poisson
(VP) simulations. Arising from collisionless plasma applications, the VP system,
\begin{equation}
\frac{\partial f}{\partial t} + {\bf v} \cdot \nabla_{\bf x} f  +
\mathbf{E}({\bf x},t) \cdot \nabla_{\bf v} f = 0, \label{eq: vlasov}
\end{equation}
and
\begin{equation}
\mathbf{E}(\mathbf{x},t)=-\nabla_{\bf x}\phi(\mathbf{x},t),\quad
-\Delta_{\bf
x}\phi(\mathbf{x},t)=\rho(\mathbf{x},t)-1,\label{eq: poisson}
\end{equation}
describes the temporal evolution of the particle distribution function in six dimensional phase space. $f( {\bf x},{\bf v},t)$ is
probability distribution function which describes the probability of
finding a particle with velocity $\bf{v}$ at position $\bf{x}$ at
time $t$, $\bf{E}$ is the electric field, and $\phi$ is the
self-consistent electrostatic potential. The probability
distribution function couples to the long range fields via the
charge density, $\rho(t,x) = \int_{\mathbb{R}^3} f(x,v,t)dv$,
where we take the limit of uniformly distributed infinitely massive
ions in the background. In this paper, we consider the VP system with
1-D in ${\bf x}$ and 1-D in ${\bf v}$. 

Many different approaches have been proposed for the VP simulations. 
There are the Lagrangian particle-in-cell (PIC) methods, which have been very popular
in practical high dimensional simulations due to its relatively low computational cost 
\cite{friedman1991multi, jacobs2009implicit, heikkinen2008full}. However, the Lagrangian particle approach is known to suffer
the statistical noise which is of order $1/\sqrt{N}$, where $N$ is the number of particles in a simulation. 
There are very high order Eulerian finite difference \cite{zhou2001numerical}, finite volume \cite{banks2010new}, finite element discontinuous Galerkin method 
\cite{heath2010discontinuous, cheng2011positivity}. Eulerian methods can be designed to be highly accurate in both space and in time, thus being able to resolve complicated solution structures in a more efficient manner by using a set of relatively coarse numerical mesh. However, they are subject to CFL time step restrictions.
There are the dimensional split SL approach originally proposed in
\cite{cheng}, and further developed in the finite volume \cite{FilbetSB, sonnendruecker, begue1999two, besse2003semi, crouseilles2010conservative},
finite difference \cite{carrillo2007nim, Qiu_Christlieb, Qiu_Shu2}, finite element discontinuous Galerkin framework \cite{Qiu_Shu_DG,rossmanith2011positivity} and a hybrid finite different-finite element framework \cite{Guo_Qiu}. The semi-Lagrangian framework allows for extra large numerical time steps compared with Eulerian approach, leading to some savings in computational cost. The dimensional splitting allows for a very simple implementation procedure for tracing characteristics; however it causes a second order operator splitting error in time. For convergence estimate for the semi-Lagrangian methods for the VP simulations, we refer to \cite{charles2013enhanced}. If the splitting is not performed properly, numerically instabilities are observed \cite{huot2003instability}.
In \cite{Guo_Qiu2}, an integral deferred correction method is proposed for the dimensional split SL approach to reduce the splitting error. 

In this paper, we proposes a high order truly multi-dimensional SL finite difference approach for solving the VP system. The `truly multi-dimensional' means that no operator splitting is involved. The difficulty is the tracing of characteristics with high order temporal accuracy in a time step. Especially the evolution of characteristics is due to the electric field induced by the unknown particle distribution function $f$ in the Vlasov equation \eqref{eq: vlasov}. A high order two-stage multi-derivative predictor-corrector algorithm is proposed to build up a high order characteristic-tracing algorithm based on lower order ones, with the help of moment equations of the VP system. A high order WENO interpolation is proposed to recover information among grid points. The proposed algorithm is of high order accuracy in both space and in time. However, there is no mass conservation. We discuss such issues as well as the computational cost of the proposed algorithm.

The paper is organized as follows. Section~\ref{sec2} describes the high order SL finite difference approach without operator splitting. High order way of tracing characteristics are proposed and analyzed. Issues related to computational cost and mass conservation are discussed.
Section~\ref{sec4} presents numerical simulation results. Finally, the conclusion is given in Section~\ref{sec5}.


\section{Truly multi-dimensional SL algorithm.}
\label{sec2}
\setcounter{equation}{0}
\setcounter{figure}{0}
\setcounter{table}{0}

\subsection{Algorithm framework}
\label{sec2.1}
Our goal is to design a high order SL finite difference scheme for the VP system without operator splitting. 
Consider the VP system \eqref{eq: vlasov} with 1-D in $x$ and 1-D $v$. The 2-D ${x-v}$ plane is discretized into uniformly spaces rectangular meshes, 
\[
x_{\frac12} < x_{1+\frac12} < \cdots < x_{i+\frac12}<\cdots < x_{n_x+ \frac12},
\]
\[
v_{\frac12} < v_{1+\frac12} < \cdots < v_{j+\frac12}<\cdots < v_{n_v+ \frac12}.
\]
The center of each of the rectangular cell $[x_{i-\frac12}, x_{i+\frac12}] \times [v_{j-\frac12}, v_{j+\frac12}]$ is denoted as $(x_i, v_j)$. 
We consider evolving the numerical solution $f^n_{i, j}$, $i = 1, \cdots n_x, j =1, \cdots, n_v$, where $f^n_{i, j}$
denotes the numerical solution at $(x_i, v_j)$ at the time level $t^n$.
The proposed SL algorithm in updating the solution $f^{n+1}_{ij}$ consists of the following steps. 
\begin{enumerate}
\item Characteristics are traced backward in time to $t^n$. Let the foot of the characteristic at the time level $t^n$ emanating from $(x_i, v_j)$ at $t^{n+1}$  be denoted as $(x^{\star}_i, v^{\star}_j)$. It is approximated by numerically solving the following final value problem
\begin{equation}
\label{eq: char}
 \left \{
 \begin{array}{l}
\frac{d{x}(t)}{dt} = {v(t)}, \\[2mm]
\frac{d{v(t)}}{dt} = {E(x(t), t)},\\[2mm]
x(t^{n+1}) = x_i, \\[2mm]
v(t^{n+1}) = v_j.
\end{array}
\right.
\end{equation}
Here, we remark that solving \eqref{eq: char} with high order temporal accuracy is non-trivial. Especially,  
the electric field ${\bf E}$ depends on the unknown function $f$ via the Poisson's equation \eqref{eq: poisson} in a global rather than local fashion. Moreover, being a final value problem, the electrical field $E$ is known initially only at the time step $t^n$. 
In Section~\ref{sec2.2}, we discuss the proposed high order (up to third order) way of tracing characteristics in time.
\item The solution is updated as
\beq
\label{eq: update}
f^{n+1}_{i, j} = f(x^{n, (l)}_i, v^{n, (l)}_j, t^n) \approx f(x^{\star}_i, v^{\star}_j, t^n).
\eeq
We propose to recover $f(x^{n, (l)}_i, v^{n, (l)}_j, t^n)$ by a high order (up to sixth order) WENO interpolation from $f^n_{i, j}$, $i = 1, \cdots n_x, j =1, \cdots, n_v$. The procedures are discussed in
Section~\ref{sec2.3}.
\end{enumerate}

\subsection{Tracing characteristics with high order temporal accuracy}
\label{sec2.2}

It is numerically challenging to design a one-step method to locate the foot of characteristics with high order accuracy in time. 
The electric field $E$ is not explicitly unknown; it is induced by the unknown function $f$ via the Poisson's equation \eqref{eq: poisson}.
Since it is difficult to evaluate the electric field $E$ (r.h.s. of equation \eqref{eq: char}) for some intermedia time stages between $[t^n, t^{n+1}]$, Runge-Kutta methods can't be used directly. 

Below we describe our proposed predictor-corrector procedure for locating the foot of characteristics. We will first describe a first order scheme in tracing characteristics; the second scheme is built upon the first order prediction; and the proposed third order scheme is built upon the second order prediction. 
In our notations, the superscript $^n$ denotes the time level, the subscript $i$ and $j$ denote the location $x_i$ and $v_j$ in $x$ and $v$ directions respectively, the superscript $^{(l)}$ denotes the formal order of approximation. For example, in equation \eqref{eq: x_v_1} below,
$x^{n, (1)}_{i}$ (or $v^{n, (1)}_{j}$) approximates $x_i^\star$ (or $v_j^\star$) with first order,
and $E^n_i = E(x_i, t^n)$. $\frac{d}{dt} = \frac{\partial}{\partial t} +  \frac{\partial x}{\partial t}\frac{\partial}{\partial x}$ denotes the material derivatives along characteristics. 
The order of approximation we mentioned in this subsection is for temporal accuracy. 
We propose to use a spectrally accurate fast Fourier transform (FFT) in solving the Poisson's equation \eqref{eq: poisson}, whose r.h.s. function $\rho(x, t) = \int f(x, v, t)dv$ is evaluated by a mid-point rule numercally. The mid point rule is of spectral accuracy given the function being integrated is either periodic or compactly supported \cite{boyd2001caf}. 

\bigskip
\noindent
\underline{\em \bf First order scheme.} We let
\beq
\label{eq: x_v_1}
x^{n, (1)}_i = x_i - v_j \Delta t; \quad v^{n, (1)}_j = v_j - E^n_i \Delta t,  
\eeq
which are first order approximations to $x_i^\star$ and $v_j^\star$, see Proposition~\ref{prop: order1} below.
Let
\beq
\label{eq: f_1}
f^{n+1, (1)}_{i, j} = f(x^{n, (1)}_i, v^{n, (1)}_j, t^n), 
\eeq
which is a first order in time approximation to $f^{n+1}_{i, j}$.
Note that the spatial approximation in equation \eqref{eq: f_1} (and in other similar equations in this subsection) is performed via high order WENO interpolation discussed in Section~\ref{sec2.3}.
Based on $\{f^{n+1, (1)}_{i, j}\}$, we computed 
\[
\rho^{n+1, (1)}_i, \quad E^{n+1, (1)}_i
\] 
by using a mid-point rule and FFT based on the Poisson's equation \eqref{eq: poisson}.
Note that  $\rho^{n+1, (1)}_i$ and $E^{n+1, (1)}_i$ also approximate $\rho^{n+1}_i$ and $E^{n+1}_i$ with first order temporal accuracy. 
\begin{prop}
\label{prop: order1}
$x^{n, (1)}_i$ and $v^{n, (1)}_j$ constructed in equation \eqref{eq: x_v_1} are first order approximations to $x_i^\star$ and $v_j^\star$ in time. 
\end{prop}
\noindent
{\em Proof.} 
By Taylor expansion,
\beqa
x_i^\star &=&  x_i - \frac{d x_i}{dt}(x_i, v_j, {t^{n+1}}) \Delta t +  \mathcal{O}(\Delta t^2) \nonumber\\
&=& x_i - v_j \Delta t +  \mathcal{O}(\Delta t^2) \nonumber\\
&\stackrel{\eqref{eq: x_v_1}}{=}& x^{n, (1)}_i +  \mathcal{O}(\Delta t^2), \nonumber
\eeqa
\beqa
v_j^\star  &=&  v_j - \frac{d  v_j}{dt}|_{t^{n+1}} \Delta t + \mathcal{O}(\Delta t^2) \nonumber\\
&=& v_j - E^{n+1}_i \Delta t +   \mathcal{O}(\Delta t^2) \nonumber\\
&=& v_j - (E^{n}_i +  \mathcal{O}(\Delta t)) \Delta t +   \mathcal{O}(\Delta t^2) \nonumber\\
&\stackrel{\eqref{eq: x_v_1}}{=} & v^{n, (1)}_j + \mathcal{O}(\Delta t^2). \nonumber
\eeqa
Hence $x^{n, (1)}_i$ and $v^{n, (1)}_j$ are second order approximations to $x_i^\star$ and $v_j^\star$ locally in time for a time step; the approximation is of first order in time globally. We remark that the proposed first order scheme is similar to, but different from, the standard forward Euler or backward Euler integrator. It is specially tailored to the system \eqref{eq: char}. 
$\Box$

\bigskip
\noindent
\underline{\em \bf Second order scheme.} We let
 \beq
 \label{eq: x_v_2}
  x^{n, (2)}_i = x_i - \frac12 (v_j + v_j^{n, (1)}) \Delta t,
  \quad
  v^{n, (2)}_j = v_j - \frac12 (E(x_i^{n, (1)}, t^n) + E^{n+1, (1)}_i)\Delta t,  
 \eeq
which are second order approximations to $x_i^\star$ and $v_j^\star$, see Proposition~\ref{prop: order2} below. 
Note that $E(x_i^{n, (1)}, t^n)$ in equation \eqref{eq: x_v_2} can be approximated by WENO interpolation from  $\{E^n_i\}_{i=1}^{n_x}$.
Let
$
f^{n+1, (2)}_{i, j} = f(x^{n, (2)}_i, v^{n, (2)}_j, t^n), 
$
approximating $f^{n+1}_{i, j}$ with second order in time.
Based on $\{f^{n+1, (2)}_{i, j}\}$, we computed 
$
\rho^{n+1, (2)}_i, \quad E^{n+1, (2)}_i
$ 
approximating $\rho^{n+1}_i$ and $E^{n+1}_i$ with second order temporal accuracy. 
\begin{prop}
\label{prop: order2}
$x^{n, (2)}_i$ and $v^{n, (2)}_j$ constructed in equation \eqref{eq: x_v_2} are second order approximations to $x_i^\star$ and $v_j^\star$ in time. 
\end{prop}
\noindent
{\em Proof.} 
It can be checked by Taylor expansion
\beqa
x_i^\star &=& x_i - \left(\frac{d x}{dt}(x_i, v_j, {t^{n+1}}) + \frac{d x}{dt}(x_i^\star, v_j^\star, t^{n})\right) \frac{\Delta t}{2} + \mathcal{O}(\Delta t^3) \nonumber\\
&=& x_i - \left(v_j^\star + v_j \right) \frac{\Delta t}{2} + \mathcal{O}(\Delta t^3) \nonumber\\
&\stackrel{Prop. \ref{prop: order1}}{=}& x_i - \left(v_j^{n, (1)} +  \mathcal{O}(\Delta t^2)  + v_j \right) \frac{\Delta t}{2} + \mathcal{O}(\Delta t^3) \nonumber\\
&=& x_i - \left(v_j^{n, (1)} + v_j \right) \frac{\Delta t}{2} + \mathcal{O}(\Delta t^3) \nonumber\\
&\stackrel{\eqref{eq: x_v_2}}{=}& x^{n, (2)}_i +  \mathcal{O}(\Delta t^3). \nonumber
\eeqa
Similarly,
\beqa
v_j^\star  &=& v_j - \left(E^{n+1}_i  + E(x_i^\star, t^n)\right) \frac{\Delta t}{2} +   \mathcal{O}(\Delta t^3) \nonumber\\
&\stackrel{Prop. \ref{prop: order1}}{=}& v_j - \left(E^{n+1, (1)}_i + E(x^{n, (1)}_i, t^n)  + \mathcal{O}(\Delta t^2)\right)\frac{\Delta t}2 +   \mathcal{O}(\Delta t^3) \nonumber\\
&\stackrel{\eqref{eq: x_v_2}}{=} & v^{n, (2)}_j + \mathcal{O}(\Delta t^3). \nonumber
\eeqa
Hence $x^{n, (2)}_i$ and $v^{n, (2)}_j$ are third order approximations to $x_i^\star$ and $v_j^\star$ locally in time for a time step; the approximation is of second order in time globally. Again the proposed second order scheme tailored to the system \eqref{eq: char} is similar to, but slightly different from, the second order Runge-Kutta integrator based on the trapezoid rule. 
$\Box$

\bigskip
\noindent
\underline{\em \bf Third order scheme.} We let
\beq
\label{eq: x_3}
  x^{n, (3)}_i = x_i - v_j \Delta t + \frac{\Delta t^2}2 (\frac23 E^{n+1, (2)}_i + \frac13 E(x_i^{n, (2)}, t^n)),
\eeq
  \begin{eqnarray}
  \label{eq: v_3}
  v^{n, (3)}_j &=& v_j - E^{n+1, (2)}_i \Delta t +  \frac{\Delta t^2}2 \left(
  \frac23 (\frac{d}{dt}E(x_i, t^{n+1}))^{(2)} + \frac13 \frac{d}{dt}E(x_i^{n, (2)}, t^n)
  \right),
 \end{eqnarray}
which are third order approximations to $x_i^\star$ and $v_j^\star$, see Proposition~\ref{prop: order3} below. 
Note that $\frac{d}{dt}E$ terms on the r.h.s. of equation \eqref{eq: v_3} will be obtained by using the macro-equations described below. 
Let
$
f^{n+1, (3)}_{i, j} = f(x^{n, (3)}_i, v^{n, (3)}_j, t^n), 
$
approximating $f^{n+1}_{i, j}$ with third order in time.
Based on $\{f^{n+1, (3)}_{i, j}\}$, we computed 
$
\rho^{n+1, (3)}_i, \quad E^{n+1, (3)}_i
$ 
approximating $\rho^{n+1}_i$ and $E^{n+1}_i$ with third order temporal accuracy. 
\begin{rem}We note that the mechanism to build this third order scheme is different from Runge-Kutta methods where intermedia stage solutions are constructed. It has some similarity in spirit to the Taylor-series (Lax-Wendroff type) method, where higher order time derivatives are recursively transformed into spatial derivatives. The difference with the Lax-Wendroff type time integration is: Lax-Wendroff method only uses spatial derivatives at one time level, while the proposed method used the spatial derivatives (or its high order approximations) at both $t^n$ and $t^{n+1}$ via a predictor-corrector procedure. In a sense, the proposed method is a two-stage multi-derivative method. 
\end{rem}

With $ \frac{\partial E}{\partial x} = \rho-1$ from the Poisson's equation \eqref{eq: poisson}, to compute the Lagrangian time derivative along characteristics $\frac{d}{dt}E = \frac{\partial}{\partial t} + v \frac{\partial}{\partial x}$, we only need to numerically approximate $\frac{\partial E}{\partial t}$. Notice that if we integrate the Vlasov equation \eqref{eq: vlasov} over $v$, we have
\beq
\label{eq: moment0}
\rho_t + J_x = 0,
\eeq
where $\rho(x, t)$ is the charge density and $J(x, t) = \int f v dv$ is the current density. With the Poisson's equation \eqref{eq: poisson}, and from eq.~\eqref{eq: moment0}, we have
$
\frac{\partial}{\partial x} (E_t +  J) =0, 
$
that is $E_t + J$ is independent of the spatial variable $x$. Thus
\[
E_t +  J = \frac1L \int (E_t + J(x, t)) dx = \frac1L \int J(x, t) dx,
\]
the last equality above is due to the periodic boundary condition of the problem. 
It can be shown, by multiplying the Vlasov equation \eqref{eq: vlasov} by $v$ and performing integration in both $x$- and $v$- directions,  that 
\[
\frac{\partial}{\partial t} \int J(x, t) dx = 0,
\]
therefore
\beq
\frac{\partial}{\partial t} E(x, t) +  J  =  \frac1L \int j(x, t=0) dx \doteq \bar{J^0}, \nonumber
\eeq
where $\bar{\cdot}$ denotes one's spatial average.
Hence,
\beq
\label{eq: dt_E}
\frac{d}{dt} E = (\frac{\partial}{\partial t} + v \frac{\partial}{\partial x}) E =  \bar{J^0} - J(x, t) + v (\rho-1).
\eeq
Specifically, in equation \eqref{eq: v_3}
\beqa
(\frac{d}{dt}E(x_i, t^{n+1}))^{(2)} &=&  \bar{J^0} - J^{n+1, (2)}_i  + v_j (\rho^{n+1, (2)}_i-1),  \nonumber\\
\frac{d}{dt}E(x_i^{n, (2)}, t^n) &=&  \bar {J^0} - J (x_i^{n, (2)}, t^n)  + v_j^{n, (2)} (\rho(x_i^{n, (2)}, t^n)-1). \nonumber
\eeqa
Note that $J^{n+1, (2)}_i$ and $J^{n}_i$ can be evaluated by mid-point rule from $\{f^{n+1, (2)}_{i, j} \}$ and $\{f^{n}_{i, j} \}$ respectively with spectral accuracy in space;
while $J (x_i^{n, (2)}, t^n)$ can be numerically approximated by WENO interpolation from $J^n_i$.

\begin{prop}
\label{prop: order3}
$x^{n, (3)}_i$ and $v^{n, (3)}_j$ constructed in equation \eqref{eq: x_3}-\eqref{eq: v_3} are third order approximations to $x_i^\star$ and $v_j^\star$ in time. 
\end{prop}
\noindent
{\em Proof.} 
It can be checked by Taylor expansion
\beqa
x_i^\star &=& x_i - \frac{d x}{dt}(x_i, v_j, {t^{n+1}}) {\Delta t} + \left(\frac23 \frac{d^2 x_i}{dt^2}(x_i, v_j, {t^{n+1}}) + \frac13 \frac{d^2 x_i}{dt^2}(x_i^\star, v_j^\star, {t^{n}})\right) \frac{\Delta t^2}{2}
+\mathcal{O}(\Delta t^4) \nonumber\\
&=& x_i -  v_j  {\Delta t} +
\left(\frac23 E^{n+1}_i + \frac13 E(x^\star_i, t^n)\right) \frac{\Delta t^2}{2} +
 \mathcal{O}(\Delta t^4) \nonumber\\
&\stackrel{Prop. \ref{prop: order2}}{=}& x_i -  v_j {\Delta t} + 
 \left(\frac23 E^{n+1, (2)}_i + \frac13 E(x_i^{n, (2)}, t^n) + \mathcal{O}(\Delta t^3)\right) \frac{\Delta t^2}{2}
+\mathcal{O}(\Delta t^4) \nonumber\\
&\stackrel{\eqref{eq: x_3}}{=}& x^{n, (3)}_i +  \mathcal{O}(\Delta t^4). \nonumber
\eeqa
Similarly,
\beqa
v_j^\star  =&& v_j  - E^{n+1}_i {\Delta t} + \left(\frac23 \frac{d E}{dt}(x_i, {t^{n+1}}) + \frac13 \frac{dE}{dt}(x_i^\star, {t^{n}})\right) \frac{\Delta t^2}{2}
+\mathcal{O}(\Delta t^4) \nonumber\\
\stackrel{Prop. \ref{prop: order2}}{=}&& v_j - (E^{n+1, (2)}_i + \mathcal{O}(\Delta t^3)) {\Delta t} \nonumber\\
&&
+ 
\left(\frac23 (\frac{d E}{dt}(x_i, {t^{n+1}}))^{(2)} + \frac13 \frac{dE}{dt}(x_i^{n, (2)}, {t^{n}})+ \mathcal{O}(\Delta t^3) \right) \frac{\Delta t^2}{2}
+   \mathcal{O}(\Delta t^4) \nonumber\\
\stackrel{\eqref{eq: v_3}}{=} && v^{n, (3)}_j + \mathcal{O}(\Delta t^4). \nonumber
\eeqa
Hence $x^{n, (3)}_i$ and $v^{n, (3)}_j$ are fourth order approximations to $x_i^\star$ and $v_j^\star$ locally in time for a time step; the approximation is of third order in time globally.  
$\Box$

\bigskip
\noindent
\underline{\em \bf Higher order extensions.} The procedures proposed above for locating the foot of characteristics can be extended to schemes with higher order temporal accuracy by using higher order version of Taylor expansion, e.g. as in equation~\eqref{eq: x_3} \eqref{eq: v_3}. As higher order material derivatives, e.g. $\frac{d^2}{dt^2}E$, are involved, a set of macro-equations from the Vlasov equation are needed. Specifically, we propose to multiply the Vlasov equation \eqref{eq: vlasov} by $v^k$, integrate over $v$ and obtain
\[
\frac{\partial }{\partial t} M_k + \frac{\partial }{\partial x} M_{k+1} - k E M_{k-1} = 0,
\]
where $M_k(x, t) = \int f(x, v, t) v^k dv$. Especially, $M_0 = \rho(x, t)$ is the charge density and $M_1 = J(x, t)$ is the current density. When $k=0$, we have equation \eqref{eq: moment0}; When $k=1$, we have
\beq
\label{eq: moment1}
\frac{\partial }{\partial t} J + \frac{\partial }{\partial x} M_2 - E \rho = 0.
\eeq
With these, we have
\beqa
   \frac{d^2E}{dt^2} &\stackrel{\eqref{eq: dt_E}}{=}& (\frac{\partial }{\partial t} + v \frac{\partial }{\partial x}) (\bar{J^0} - J(x, t) + v (\rho-1)) \nonumber\\
   \label{eq: dEdt2}
   &\stackrel{\eqref{eq: moment1}}{=}&  v^2 \pad{\rho}{x} + \pad{M_2}{x} -2v\pad{J}{x} - E,
\eeqa
where spatial derivative terms can be evaluated by high order WENO interpolations or reconstructions. 



\subsection{High order WENO interpolations.}
\label{sec2.3}

In this subsection, we discuss the procedures in spatial interpolation to recover information among grid points, e.g. to update numerical solution by equation \eqref{eq: update}, and in spatial reconstruction to recover function derivatives at grid points, e.g. in computing spatial derivatives in equation ~\eqref{eq: dEdt2}. There have been a variety of interpolation choices, such as the piecewise parabolic method (PPM) \cite{colella1984piecewise}, spline interpolation \cite{crouseilles2007hermite}, cubic interpolation propagation (CIP) \cite{yabe2001cip}, ENO/WENO interpolation \cite{carrillo2007nim, Qiu_Shu2}. In our work we adapt the WENO interpolations. 

\bigskip
\noindent
\underline{\em \bf WENO interpolations.} High order accuracy is achieved by using several points in the neighborhood: the number of points used in the interpolation 
determines the order of interpolation. WENO \cite{Shu_book, carrillo2007nim, Qiu_Shu2}, short for `weighted essentially non-oscillatory',  is a well-developed adaptive procedure to 
overcome Gibbs phenomenon, when the solution is under-resolved or contains discontinuity. Specifically, when the solution is smooth the WENO interpolation recovers the linear interpolation for very high order accuracy; when the solution is under-resolved, the WENO interpolation automatically assign more weights to smoother stencils. The smoothness of the stencil is measured by the divided differences of numerical solutions. Below we provide formulas for the sixth order WENO interpolations, which is what we used in our simulations. 

The sixth order WENO interpolation at a position $x\in [x_{i-1}, x_{i}]$ (or $\xi \doteq \frac{x-x_i}{\Delta x} \in [-1, 0]$) is obtained by
\[
Q(\xi) = \omega_1 P_1(\xi) + \omega_2 P_2(\xi) + \omega_3 P_3(\xi),
\]
where
\[
P_1(\xi) =   (f_{i-3}, f_{i-2}, f_{i-1}, f_i) \,
\left (
\begin{array}{llll}
0&-1/3&-1/2&-1/6\\
0 & 3/2 & 2&1/2 \\
0&-3&-5/2&-1/2\\
1&11/6&1&1/6\\
\end{array}
\right )
\,
\left (
\begin{array}{l}
1\\
\xi\\
\xi^2\\
\xi^3
\end{array}
\right ),
\]
\[
P_2(\xi) =   (f_{i-2}, f_{i-1}, f_i, f_{i+1}) \,
\left (
\begin{array}{llll}
0&1/6&0&-1/6\\
0 & -1 & 1/2&1/2 \\
1&1/2&-1&-1/2\\
0&1/3&1/2&1/6
\end{array}
\right )
\,
\left (
\begin{array}{l}
1\\
\xi\\
\xi^2\\
\xi^3
\end{array}
\right ), 
\]
\[
P_3(\xi) = (f_{i-1}, f_i, f_{i+1}, f_{i+2}) \,
\left (
\begin{array}{llll}
0&-1/3&1/2&-1/6\\
1 & -1/2 & -1&1/2 \\
0&1&1/2&-1/2\\
0&-1/6&0&1/6
\end{array}
\right )
\,
\left (
\begin{array}{l}
1\\
\xi\\
\xi^2\\
\xi^3
\end{array}
\right ). 
\]
Linear weights
\[
\gamma_1(\xi) = \frac{1}{20}(\xi-1)(\xi-2) , \quad \gamma_2(\xi) = -\frac{1}{10}(\xi+3)(\xi-2), 
\quad \gamma_3(\xi) = \frac{1}{20}(\xi+3)(\xi+2) .
\]
Nonlinear weights are chosen to be
$$
\omega_m = \frac {\tilde{\omega}_m}
{\sum_{l=1}^3 \tilde{\omega}_l},\qquad \mbox{with} \quad
 \tilde{\omega}_l = \frac {\gamma_l}{(\varepsilon + \beta_l)^2} , \quad l = 1, 2, 3,
$$
where $\epsilon=10^{-6}$,
and the smoothness indicators 
\beqa
\beta_1 =
-9\,f_{{i-3}}f_{{i-2}}+4/3\,{f_{{i-3}}}^{2}-11/3\,f_{{i-3}}f_{{i}}+10\,f_{{i-3}}
f_{{i-1}}+14\,f_{{i-2}}f_{{i}}\nonumber\\
+22\,{f_{{i-1}}}^{2}-17\,f_{{i-1}}f_{{i}}+10/3\,
{f_{{i}}}^{2}+16\,{f_{{i-2}}}^{2}-37\,f_{{i-2}}f_{{i-1}},\nonumber
\eeqa
\beqa
\beta_2 =
-7\,f_{{i-2}}f_{{i-1}}+4/3\,{f_{{i-2}}}^{2}-5/3\,f_{{i-2}}f_{{i+1}}+6\,f_{{i-2}}U_
{{i}}+6\,f_{{i-1}}f_{{i+1}}\nonumber\\
+10\,{f_{{i}}}^{2}-7\,f_{{i}}f_{{i+1}}+4/3\,{f_{{
4}}}^{2}+10\,{f_{{i-1}}}^{2}-19\,f_{{i-1}}f_{{i}},\nonumber
\eeqa
\beqa
\beta_3 = 
-17\,f_{{i-1}}f_{{i}}+10/3\,{f_{{i-1}}}^{2}-11/3\,f_{{i-1}}f_{{i+2}}+14\,f_{{i-1
}}f_{{i+1}}+10\,f_{{i}}f_{{i+2}}\nonumber\\
+16\,{f_{{i+1}}}^{2}-9\,f_{{i+1}}f_{{i+2}}+4/3\,
{f_{{i+2}}}^{2}+22\,{f_{{i}}}^{2}-37\,f_{{i}}f_{{i+1}}.\nonumber
\eeqa

\subsection{Computational cost and savings}

One of the procedures in the proposed algorithm that takes up much computational time is to trace the foot of characteristics. Assume $N = n_x = n_v$, the scheme involves solving the Poisson's equation via FFT with the cost on the order of $N log(N)$ and a high order 2-D WENO interpolation on the order of $C N^2$, where the constant $C$ is larger when the order of interpolation is higher. Since the 2-D WENO interpolation (compared with the 1-D Poisson solver) is a procedure that takes most of the computational time, we will use the number of 2-D WENO interpolations as a measurement of computational cost. 

For the first order scheme \eqref{eq: x_v_1}, there is a high order 2-D WENO interpolation involved. The proposed second order scheme \eqref{eq: x_v_2} is based on the first order prediction: two high order 2-D WENO interpolations are involved. This leads to twice the computational cost as a first order scheme. The third order scheme \eqref{eq: x_3} - \eqref{eq: v_3} is based on the second order prediction: three high order 2-D WENO interpolations are involved. We claim that proposed high order procedures are computationally efficient: the computational cost roughly grows linearly with the order of approximation. To further save some computational cost, we propose to use lower order 2-D WENO interpolation in the prediction steps. Specifically, in the third order scheme  \eqref{eq: x_3} - \eqref{eq: v_3}, we propose to use a second order 2-D WENO interpolation in the first order prediction, use a fourth order 2-D WENO interpolation in the second order prediction, and use a sixth order 2-D WENO interpolation in the final step of updating.

\subsection{Discussion on mass conservative correction and stability}
\label{sec3}
\setcounter{equation}{0}
\setcounter{figure}{0}
\setcounter{table}{0}

The proposed scheme is non mass conservative. One possible remedy is a conservative correction procedure, that allows the construction of a conservative scheme starting from a non conservative one. 
This approach was first introduced in the context if the BGK model of rarefied gas dynamics by P.~Santagati in his PhD thesis \cite{Santagati07}, and illustrated in a preprint 
\cite{Russo-Santagati-BGK-11}. Take a simple linear convection equation in one space dimension for example, the equation will take the form
\begin{equation}
\pad{f}{t} + \pad{f}{x} = 0, \quad f(x,0) = f^0(x),
\label{eq:scalar2}
\end{equation}
with periodic boundary conditions. \rf{eq:scalar2} is discretized on a spatial grid, $x_i = i\Dx$, $i=1\ldots,n$. 

Following Osher and Shu \cite{ShuOsherEfficient}, we impose that the pointwise value $f_i^n\approx f(x_i,t^n)$ satisfies the equation 
\[
f^{n+1}_i - f^n_i 
   = -\frac{\hF_\jp-\hF_\jm}{\Dx},
\] 
where the function $\hF$ is reconstructed at the edge of the cell from the point wise values of $F(x_i) = \int_{t^n}^{t^{n+1}} f(x_i, \tau)d\tau$ in the same way pointwise values of a function $u(x\pm\Dx)$ can be reconstructed from cell average $\bar{u}_i$, see \cite{Jiang_Shu} for a detailed description of the WENO reconstruction procedure. Let $(c_\ell,b_\ell)$, $\ell = 1,\ldots,s$ be the nodes and weights of an accurate quadrature formula in the interval $[0,1]$. 
To approximate $F(x_i)$, one can use a quadrature rule  
\[
F(x_i) \approx \Dt \sum_{\ell=1}^{s} b_\ell f(x_i, t^n + c_\ell \Dt),
\]
where $f(x_i, t^n + c_\ell \Dt)$ can be obtained by the characteristics tracing as well as WENO interpolation described earlier this section. Such procedure can be directly extended to two dimensional problem, including the Vlasov-Poisson procedure, where the non-conservative semi-Lagrangian method previously proposed can be used to get the solution at quadrature points.  The 2-point Gauss-Legendre quadrature formula with $b_1 = b_2 = 1$ and $c_{1,2}= \frac{1}{2}\pm\frac{1}{2\sqrt{3}}$ is found to be a good choice with good stability property. On the other hand, such conservative correction is subject to a time step constraint related to the spatial mesh size similar to that of the Eulerian approach from spatial interpolation and reconstruction procedures. As a result, the advantage of using larger time steps in a SL method is lost. To investigate and improve such stability constraint is subject to our future research.

\section{Numerical tests: the Vlasov-Poisson system}
\label{sec4}

In this section, we examine the performance of the
proposed fully multi-dimensional semi-Lagrangian method for the VP systems. Periodic
boundary condition is imposed in x-direction, while zero boundary
condition is imposed in v-direction. We recall several norms in the VP system
below, which should remain constant in time.
\begin{enumerate}
\item $L^p$ norm $1\leq p<\infty$:
\begin{equation}
\|f\|_p=\left(\int_v\int_x|f(x,v,t)|^pdxdv\right)^\frac1p.
\end{equation}
\item Energy:
\begin{equation}
\text{Energy}=\int_v\int_xf(x,v,t)v^2dxdv + \int_xE^2(x,t)dx,
\end{equation}
where $E(x,t)$ is the electric field.
\item Entropy:
\begin{equation}
\text{Entropy}=\int_v\int_xf(x,v,t)\log(f(x,v,t))dxdv.
\end{equation}
\end{enumerate}
Tracking relative deviations of these quantities numerically will be
a good measure of the quality of numerical schemes. The relative
deviation is defined to be the deviation away from the corresponding
initial value divided by the magnitude of the initial value.
In our numerical tests, we let the time step size $\Delta t = CFL \cdot \min(\Delta x/v_{max}, \Delta v/\max(E))$, where $CFL$ is specified for different runs; 
and let $v_{max} = 6$ to minimize the error from truncating the domain in $v$-direction. 
We first present the example of two stream instability. In this example, we will
demonstrate the (1) high order spatial accuracy and the high order
temporal accuracy of the proposed schemes; (2) the time evolution of overall mass and other theoretically conserved physical norms for the proposed method; (3) the performance of the proposed method in resolving solution structures.

\begin{exa} Consider two stream instability
\cite{FilbetS}, with an unstable initial distribution
function:
\begin{equation}
f(x,v,t=0)=\frac{2}{7\sqrt{2\pi}}(1+5v^2)(1+\alpha((\cos(2kx)+\cos(3kx))/1.2+
\cos(kx))\exp(-\frac{v^2}{2})
\end{equation}
with $\alpha=0.01$, $k=0.5$, the length of the domain in the x
direction is $L=\frac{2\pi}{k}$ and the background ion distribution
function is fixed, uniform and chosen so that the total net charge
density for the system is zero.
\end{exa}

\begin{figure}[htb]
\begin{center}
\includegraphics[height=2.2in,width=3.0in]{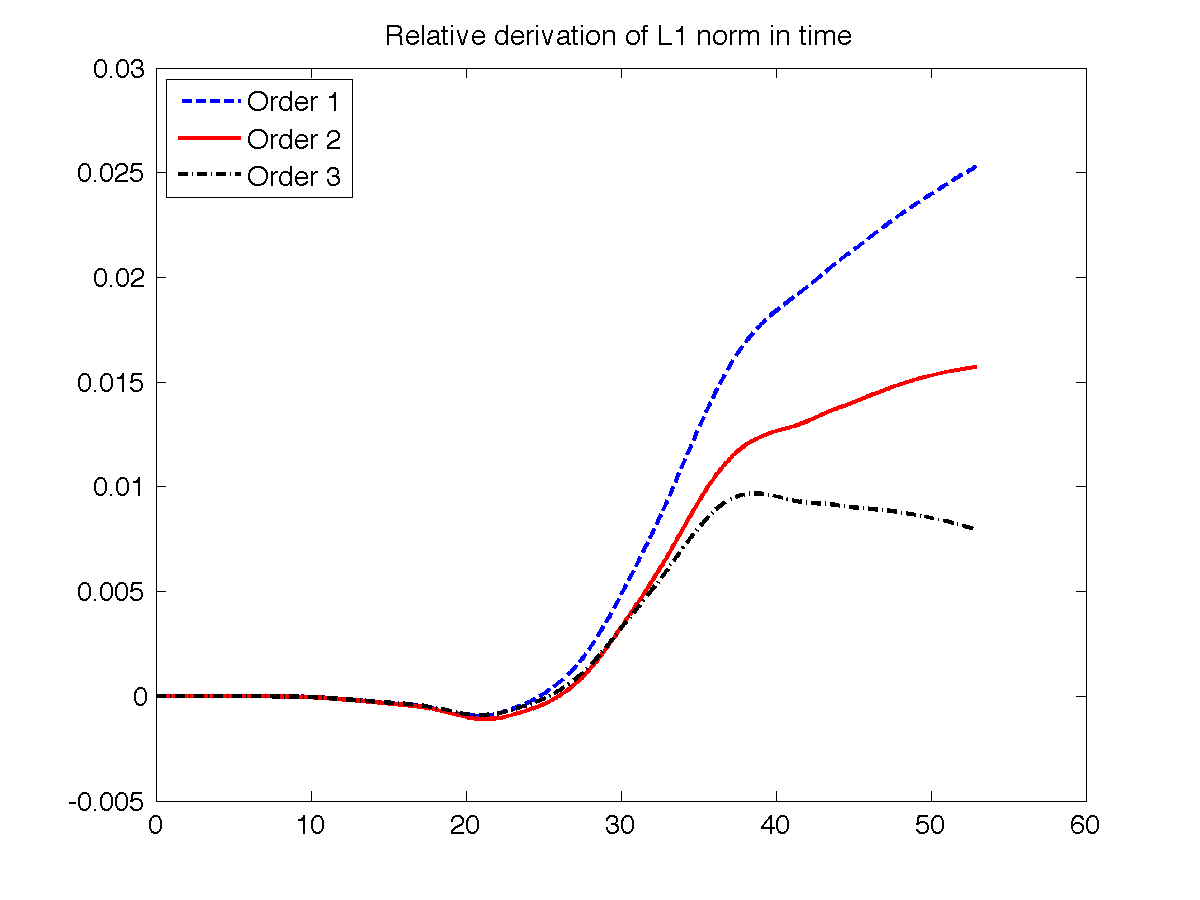}
\includegraphics[height=2.2in,width=3.0in]{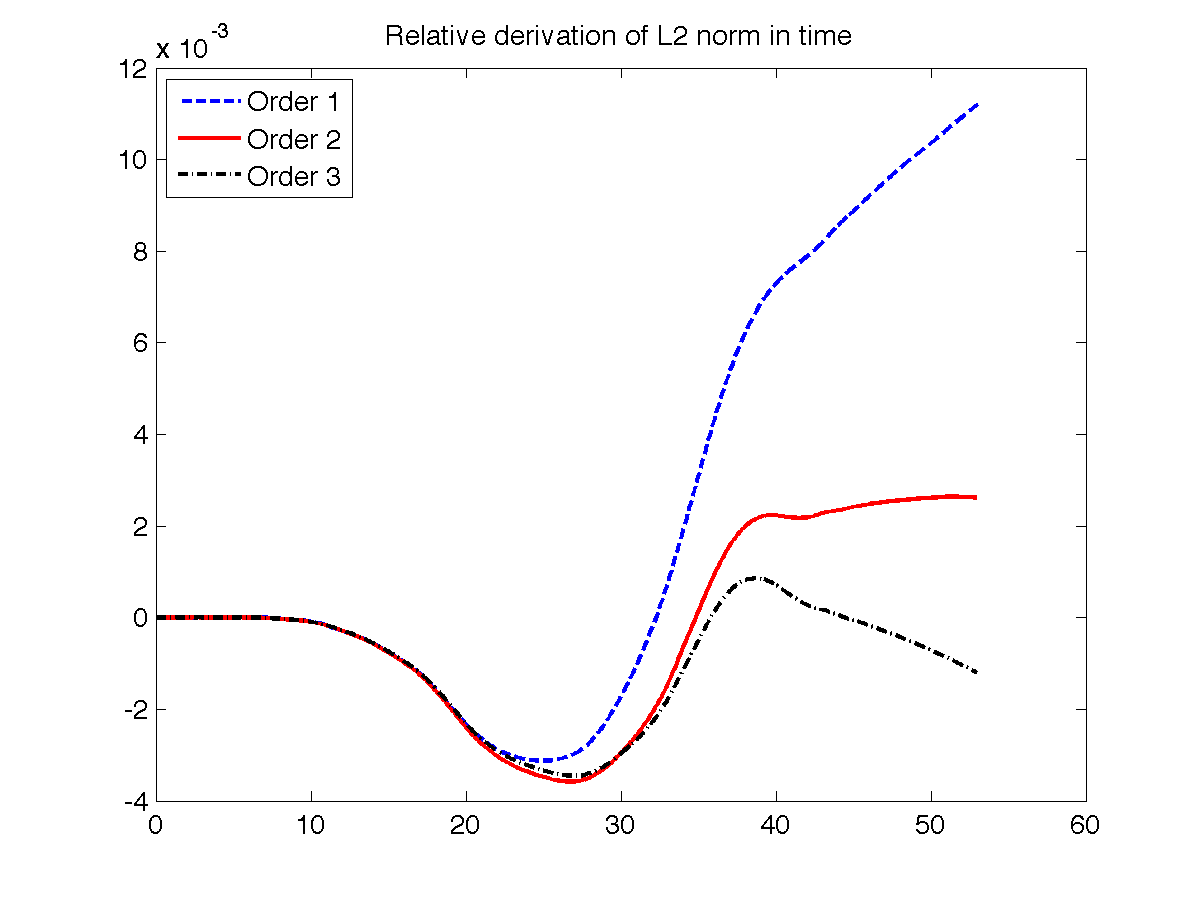}\\
\includegraphics[height=2.2in,width=3.0in]{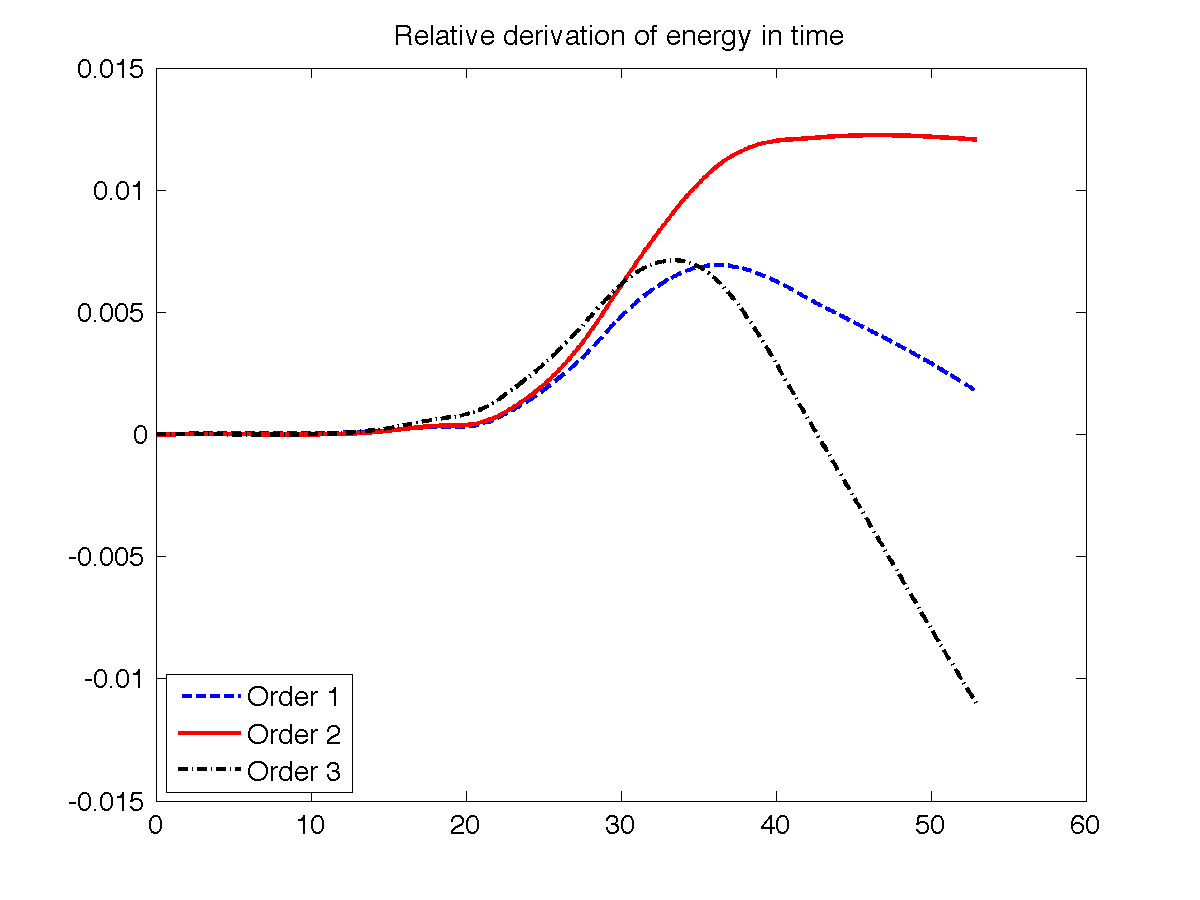}
\includegraphics[height=2.2in,width=3.0in]{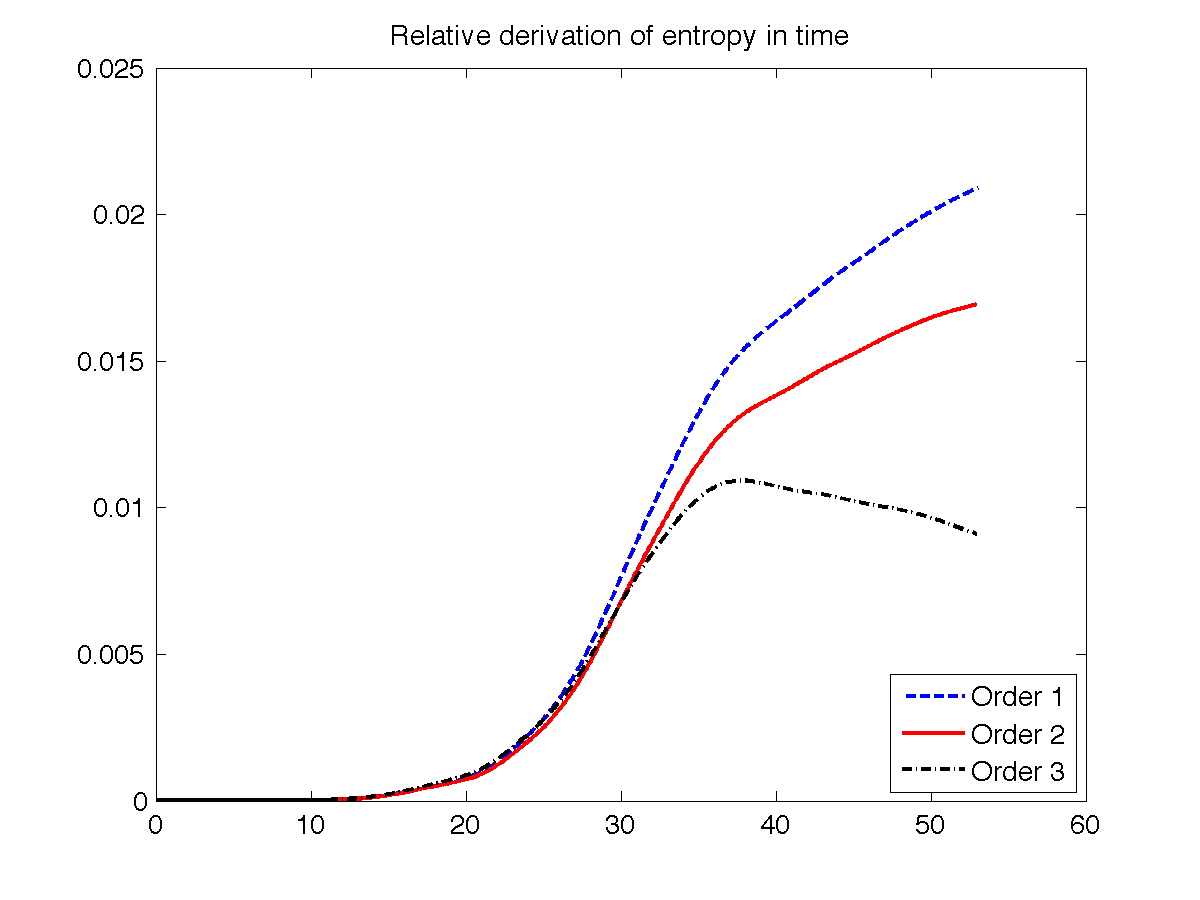}
\end{center}
\caption{Two stream instability. The SL WENO scheme with
sixth order WENO interpolation in space and various orders of temporal accuracy. 
Time evolution of the relative
deviations of discrete $L^1$ norms (upper left), $L^2$ norms,
kinetic energy norms (lower left) and entropy (lower right).}
\label{fig: 2stream_norm}
\end{figure}

\begin{figure}[htb]
\begin{center}
\includegraphics[height=3.2in,width=4.0in]{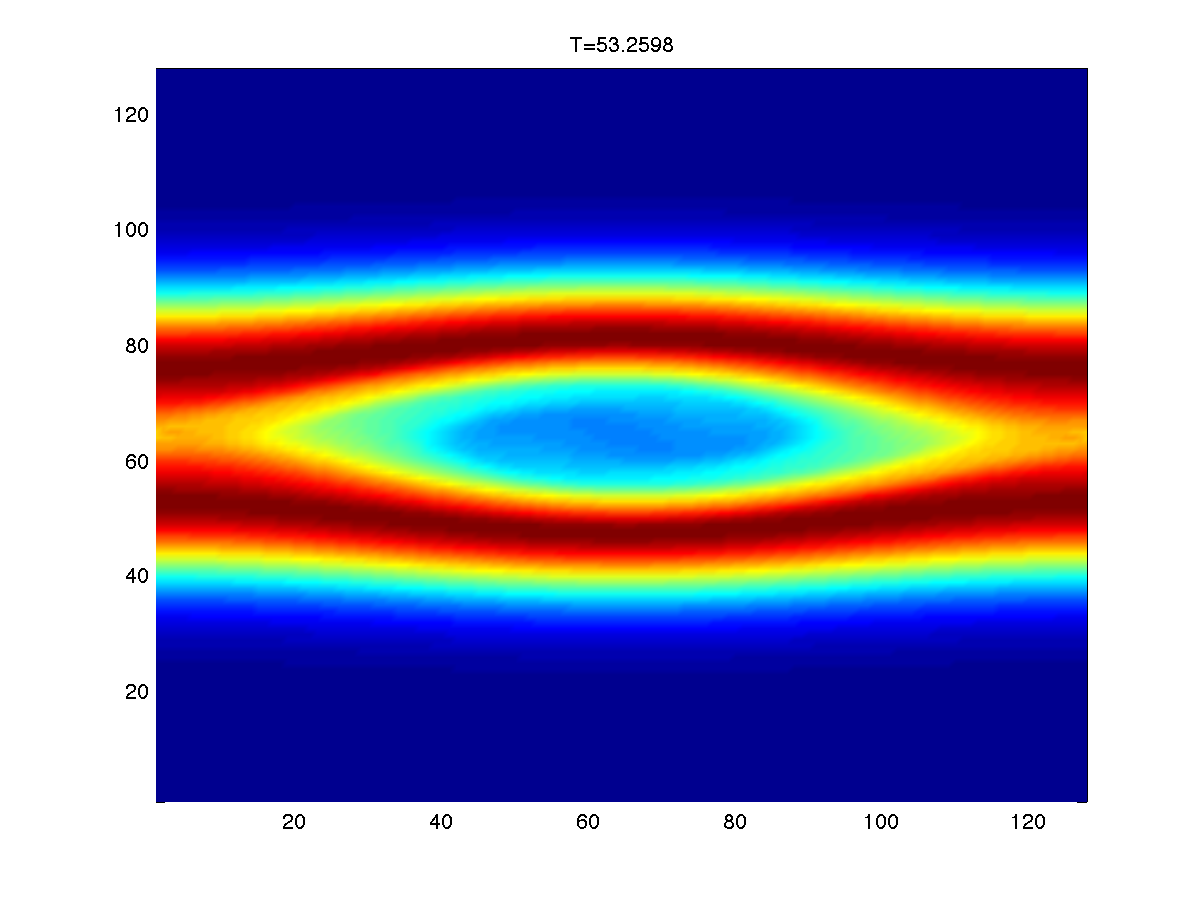}
\end{center}
\caption{Two stream instability: $T=53$. The SL WENO scheme with the sixth order WENO interpolation and a third order temporal accuracy. The spatial mesh is $128 \times 128$ and $CFL=5$.}
\label{fig: 2stream}
\end{figure}

We test both spatial and temporal convergence of the proposed truly multi-dimensional semi-Lagrangian method. 
We first test the spatial convergence by using a sequence of meshes with $n_x = n_v = \{210, 126, 90, 70\}$. 
The meshes are designed so that the coarse mesh grid coincides with part of the reference fine mesh grid ($n_x=n_v=630$).
We set $CFL=0.01$ so that the spatial error is the dominant error.
Table~\ref{tab: spa_2stream} is the spatial convergence table for the proposed schemes with sixth order WENO interpolation. The expected fifth order convergence globally in time in observed.
We then test the temporal convergence of the proposed first, second and third order schemes. 
Table~\ref{tab: tem_2stream} provides the temporal convergence rate for the scheme with the first to third order temporal accuracy. 
We use the sixth order WENO interpolation and a spatial mesh of $Nx = Nv=160$, so that the temporal error is the dominant error.
Expected first, second and third order temporal accuracy is observed. In Table~\ref{tab: tem_2stream}, the time step size is about $6$ to $10$ times that from an Eulerian method, yet highly accurate numerical results is achieved. 
To compare the performance of schemes with different temporal orders, we numerically
track the time evolution of physically conserved quantities of the system. In our runs, we let $n_x=n_v=128$, $CFL=5$. 
In Figure~\ref{fig: 2stream_norm}, the time evolution of numerical $L^1$ norm, $L^2$ norm, energy and entropy for schemes
with different orders of temporal accuracy are plot. In general, high order temporal accuracy indicates a better preservation of those physically conserved norms. The $L^1$ norm is not conserved since our scheme is neither mass conservative nor positivity preserving. In Figure \ref{fig: 2stream}, we show the contour plot of the numerical solution of the proposed SL WENO method with third order temporal accuracy at around $T=53$. The plot is comparable to our earlier work reported in \cite{Qiu_Christlieb, Qiu_Shu2}.

\begin{table}[htb]
\begin{center}
\bigskip
\begin{tabular}{|c | c c|}
\hline
\cline{1-3} $Nx \times Nv$ &$L^1$ error & order \\
\hline
{$70\times 70$} &7.01E-7 & -- \\
\hline
{$90\times 90$} &2.06E-7&4.88\\
\hline
{$126\times 126$} &3.96E-8&4.89\\
\hline
{$210\times 210$}  &3.20E-9&4.95\\
\hline
\end{tabular}
\end{center}
\caption{Order of accuracy in space for the SL WENO schemes: two stream instability.
The scheme use sixth order WENO interpolation and has a third order temporal accuracy in tracing characteristics.
$T=1$ and $CFL=0.01$.}
\label{tab: spa_2stream}
\end{table}

\begin{table}[htb]
\begin{center}
\bigskip
\begin{tabular}{|c | c c|c c|c c|}
\hline
\cline{1-7} &\multicolumn{2}{c|}{first order} &\multicolumn{2}{c|}{second order} &\multicolumn{2}{c|}{third order}  \\
\hline
\cline{1-7} $CFL$& $L^1$ error&order&  $L^1$ error&order&  $L^1$ error&order\\
\hline
6 & 1.17E-4& -- &2.40E-6 & -- & 1.13E-7&--\\
\hline
7 & 1.40E-4&1.13  & 2.80E-6 & 2.04 & 1.79E-7&3.02\\
\hline
8 & 1.63E-4& 1.16 & 3.69E-6 &2.07 & 2.69E-7&3.02\\
\hline
9 & 1.87E-4& 1.16 & 4.69E-6 &2.04 & 3.84E-7&3.03\\
\hline
10 & 2.12E-4& 1.20 & 5.84E-6 &2.08 & 5.31E-7&3.06\\
\hline
\end{tabular}
\end{center}
\caption{Order of accuracy in time for the SL WENO schemes with sixth order WENO interpolation and various orders of temporal accuracy. Two stream instability. $Nx = Nv=160$ and $T=5$.}
\label{tab: tem_2stream}
\end{table}

\begin{exa} Consider weak Landau damping for the
Vlasov-Poisson system with initial condition:
\begin{equation}
\label{landau}
f(x,v,t=0)=\frac{1}{\sqrt{2\pi}}(1+\alpha\cos(kx))\exp(-\frac{v^2}{2}),
\end{equation}
where $\alpha=0.01$. 
When the perturbation magnitude is small enough ($\alpha=0.01$), the VP system
can be approximated by linearization around the Maxwellian equilibrium
$f^0(v)=\frac{1}{\sqrt{2\pi}}e^{-\frac{v^2}{2}}$.
The analytical damping rate of electric field can be derived accordingly \cite{fried1961plasma}.
We test the numerical numerical damping rates with theoretical values. We only present
the case of $k=0.5$. The spatial computational grid has $n_x=n_v=128$ and $CFL=5$. 

For the scheme with first, second and third order accuracy in time and sixth order WENO interpolation in space, we plot the evolution of electric field in $L^2$ norm
benchmarked with theoretical values (solid black lines in the figure) in Figure \ref{fig403}. A better match with the theoretical decay rate of the electric field is observed for schemes with second and third order temporal accuracy. The time evolution of discrete $L^1$ norm, $L^2$ norm, kinetic energy and
entropy of schemes with different temporal orders are reported in Figure \ref{fig404}. $L^1$ and $L^2$ norms are better preserved by schemes with higher order temporal accuracy. 
Note that the mass is not exactly preserved. Energy and entropy are better preserved by schemes with second and third order accuracy than that with first order accuracy. 

\begin{figure}
\begin{center}
\includegraphics[height=2.2in,width=3.0in]{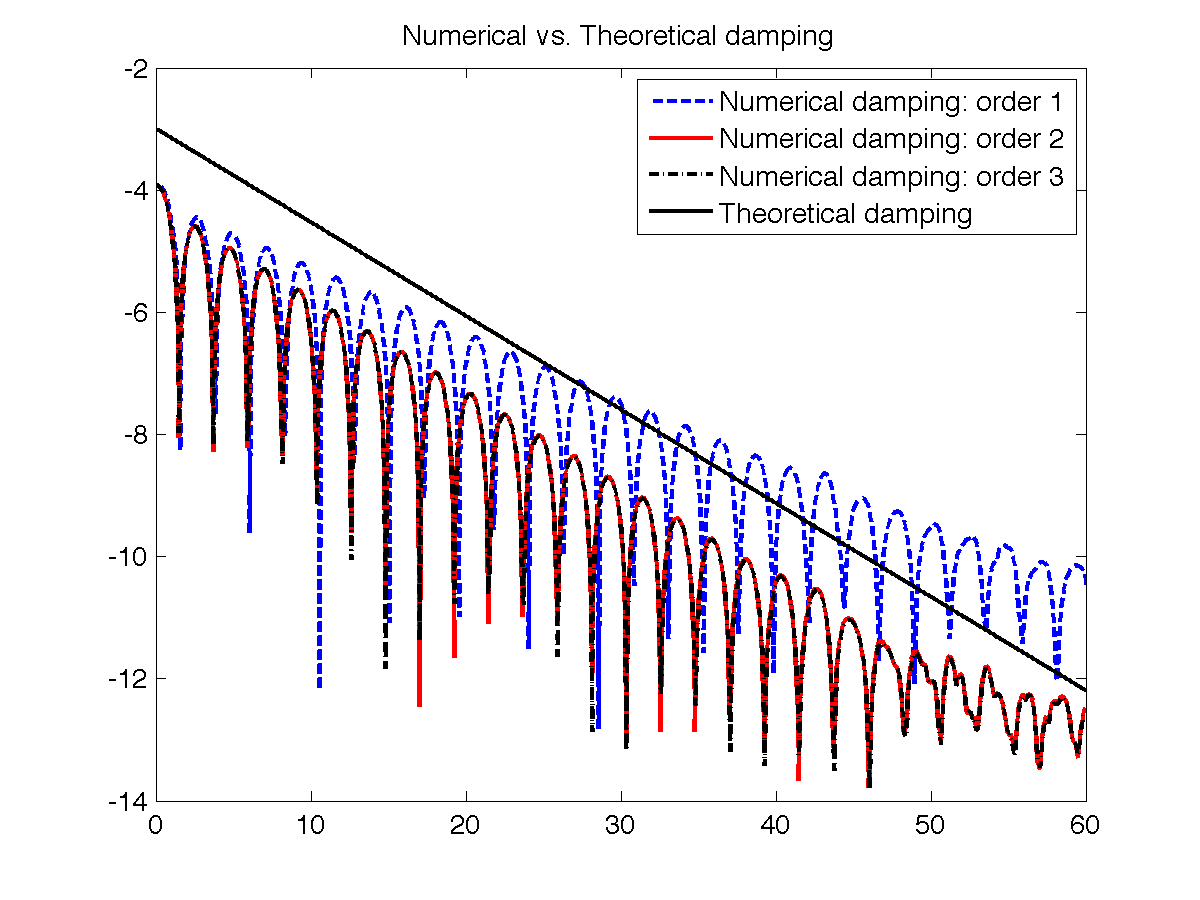}
\end{center}
\caption{Weak Landau damping. Time evolution of electric field in
$L^2$ norm.} 
\label{fig403}.
\end{figure}

\begin{figure}[htb]
\begin{center}
\includegraphics[height=2.2in,width=3.0in]{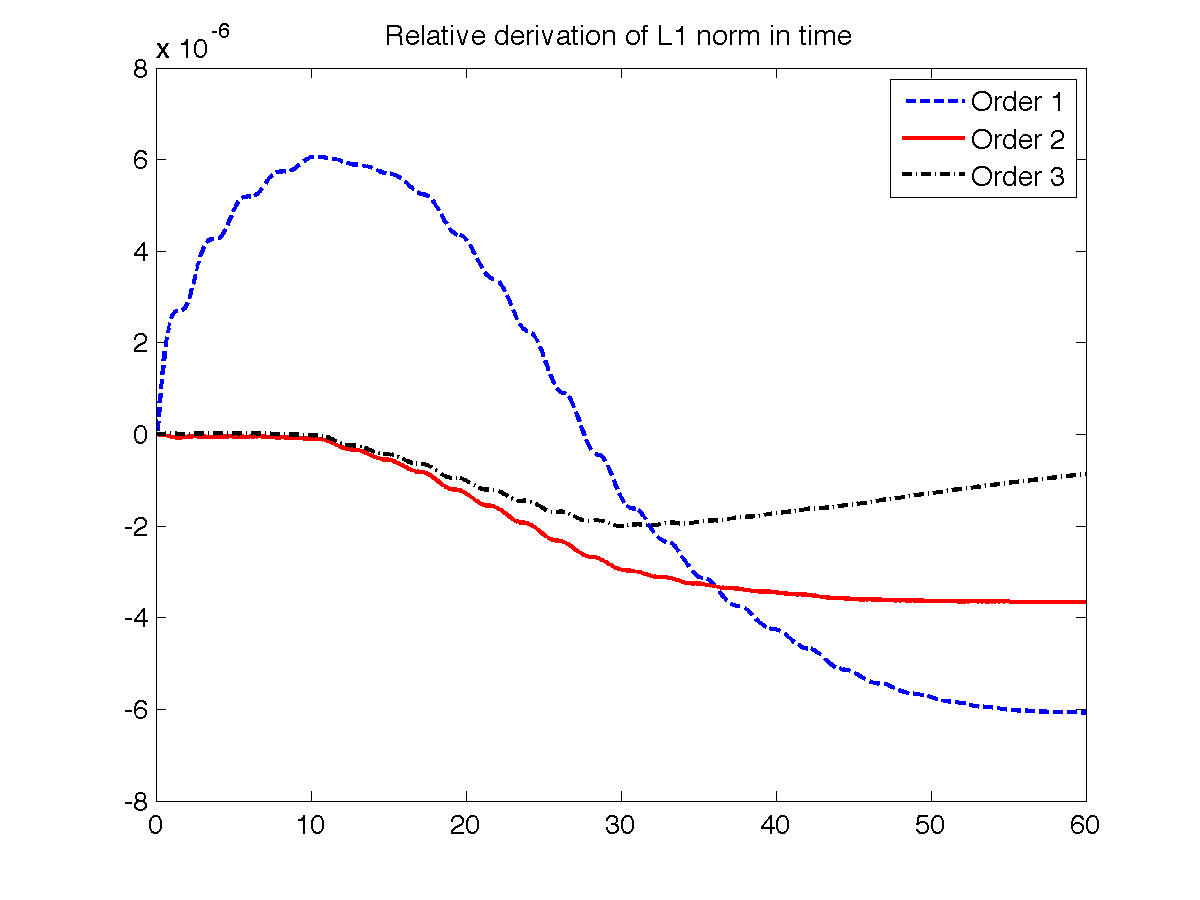}
\includegraphics[height=2.2in,width=3.0in]{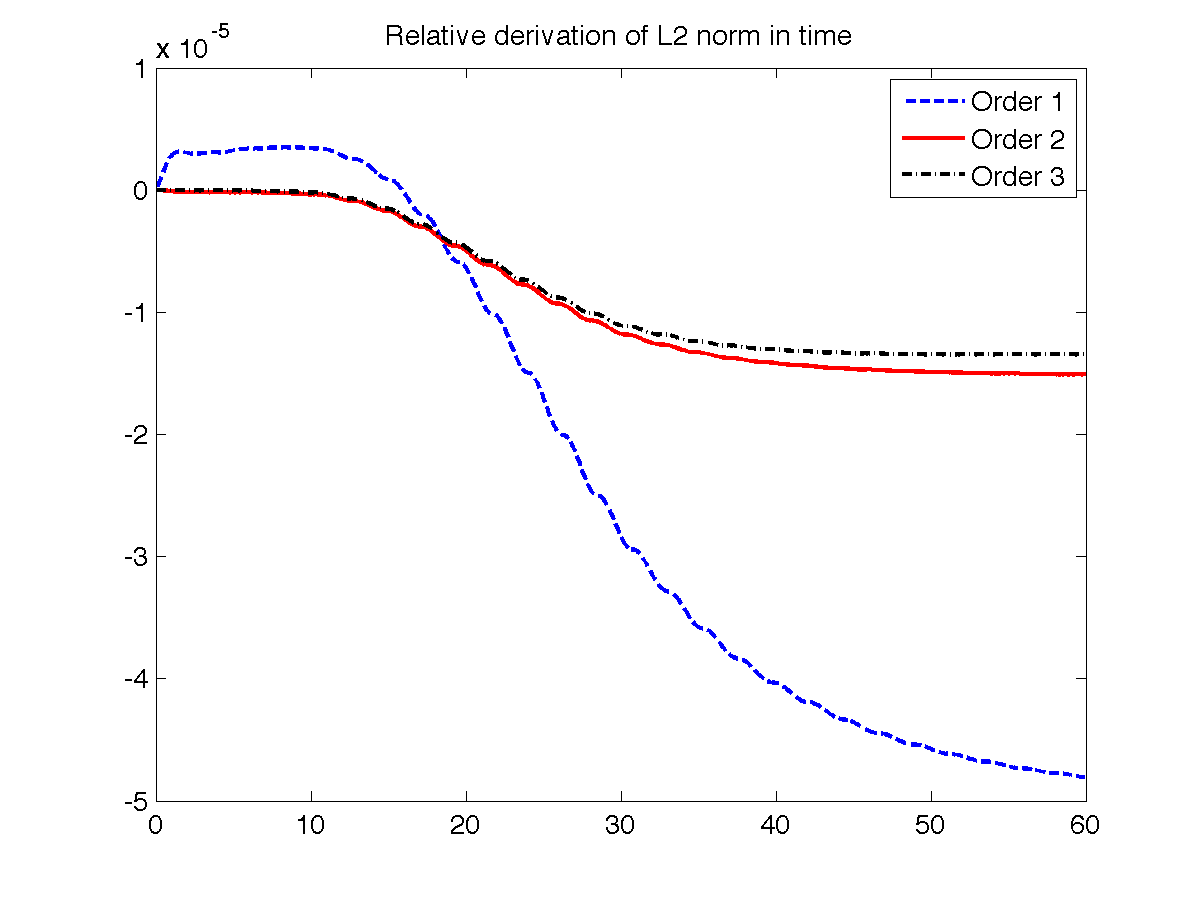}\\
\includegraphics[height=2.2in,width=3.0in]{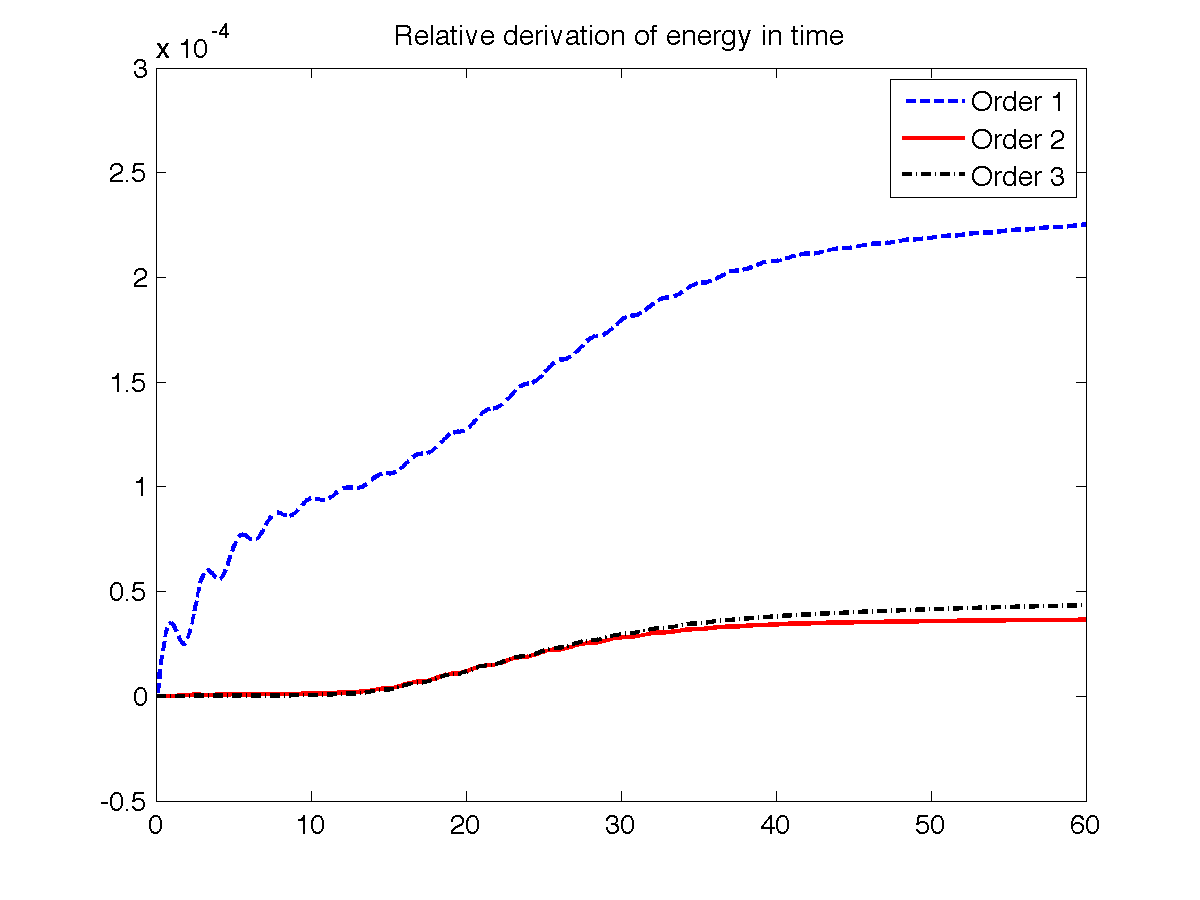}
\includegraphics[height=2.2in,width=3.0in]{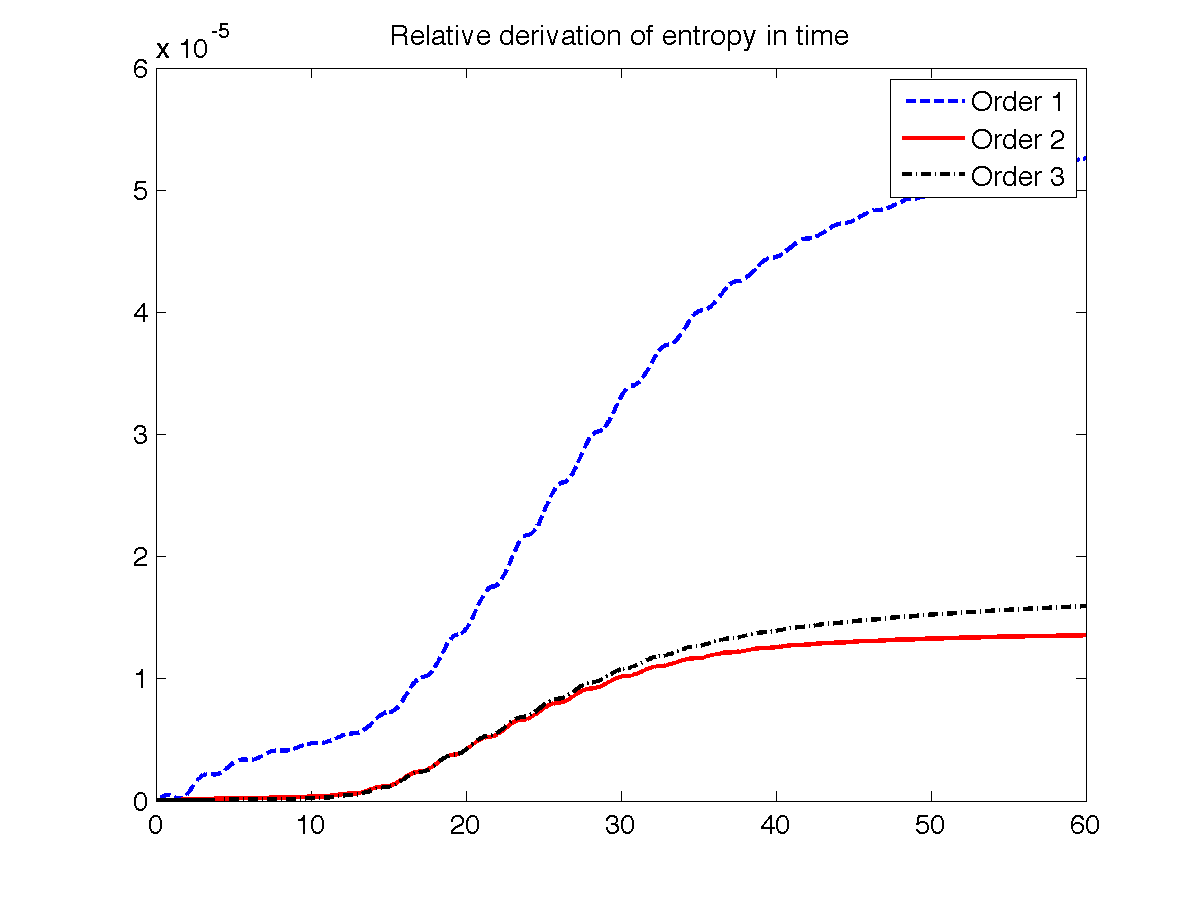}
\end{center}
\caption{Weak Landau damping. The proposed SL WENO scheme with
first, second and third order accuracy in time and sixth order WENO interpolation in space. 
Time evolution of the relative
deviations of discrete $L^1$ norms (upper left), $L^2$ norms,
kinetic energy norms (lower left) and entropy (lower right).}
\label{fig404}
\end{figure}

\end{exa}

\begin{exa} Consider strong Landau damping. The
initial condition is equation \eqref{landau}, with $\alpha=0.5$ and
$k=0.5$. 
The evolution of $L^2$ norms of electric field is provided in Figure \ref{fig405}, which is comparable to existing results in the literature, e.g. see \cite{Guo_Qiu}. 
The time evolution of discrete $L^1$ norm,
$L^2$ norm, kinetic energy and entropy are reported in Figure \ref{fig407}. 
The $L^1$ norm, as expected, is not conservative. 
Numerical solutions of the proposed scheme at different times are observed to be comparable 
to those that have been well reported in the literature, e.g. \cite{Qiu_Christlieb, Guo_Qiu} among many others. Thus we omit to present
those figures to save space. 

\begin{figure}
\begin{center}
\includegraphics[height=2.2in,width=3.0in]{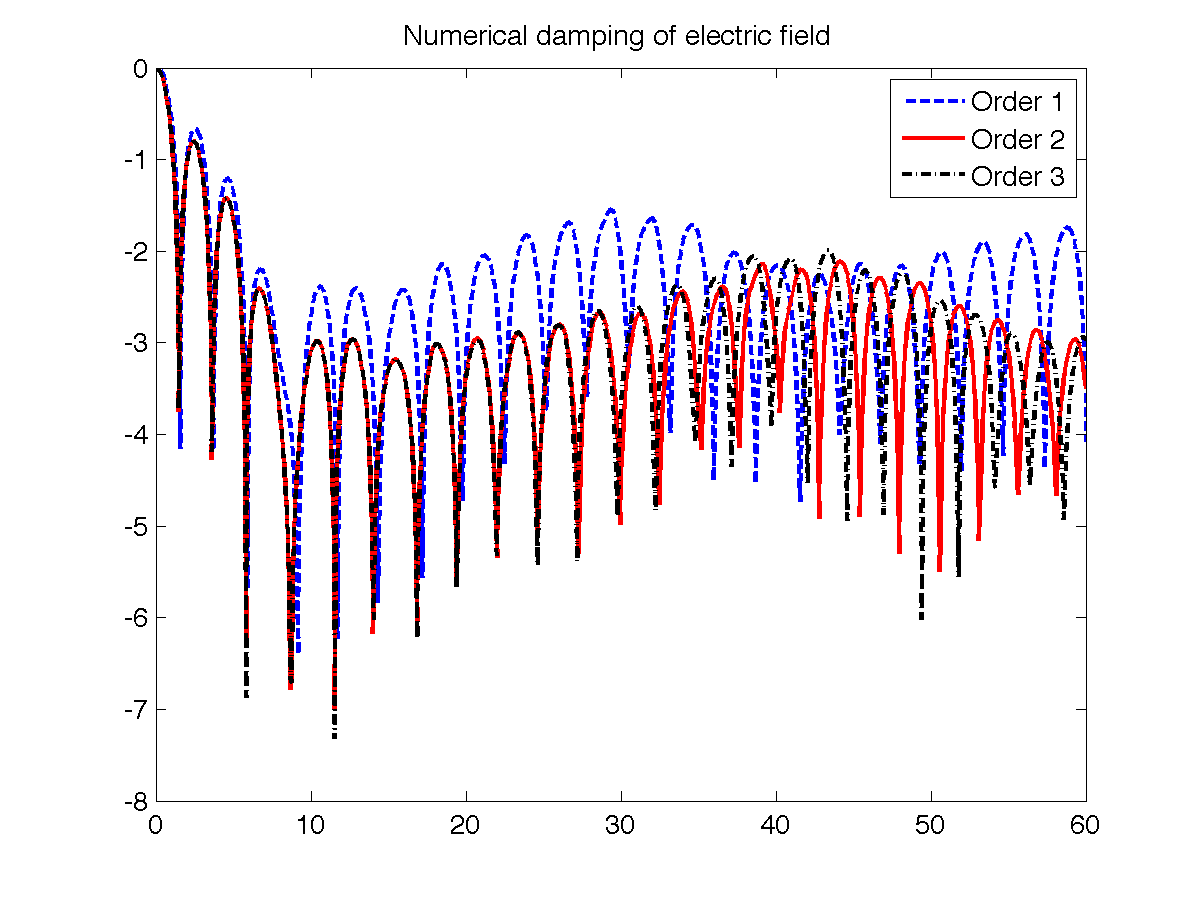}
\end{center}
\caption{Strong Landau damping. Time evolution of electric field in
$L^2$ norm.} 
\label{fig405}.
\end{figure}

\begin{figure}[htb]
\begin{center}
\includegraphics[height=2.2in,width=3.0in]{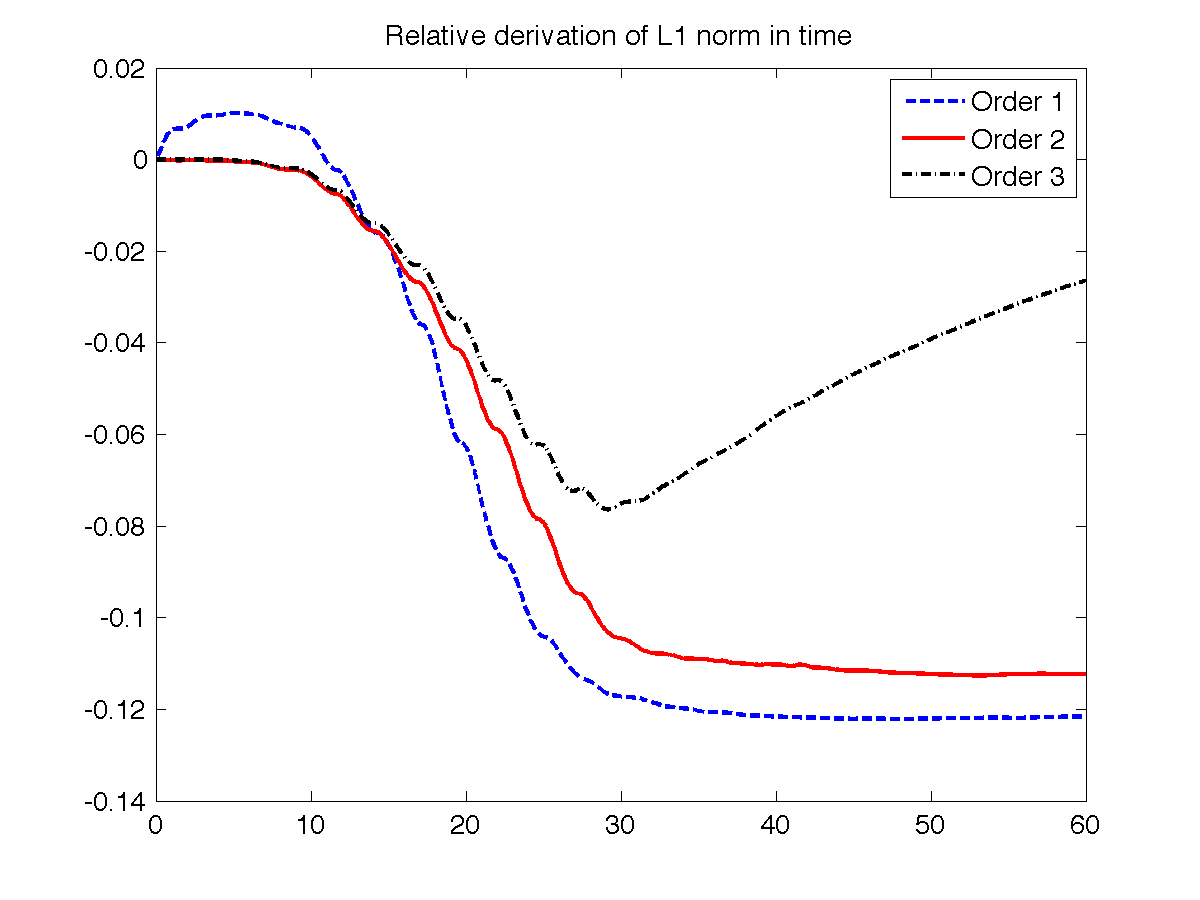}
\includegraphics[height=2.2in,width=3.0in]{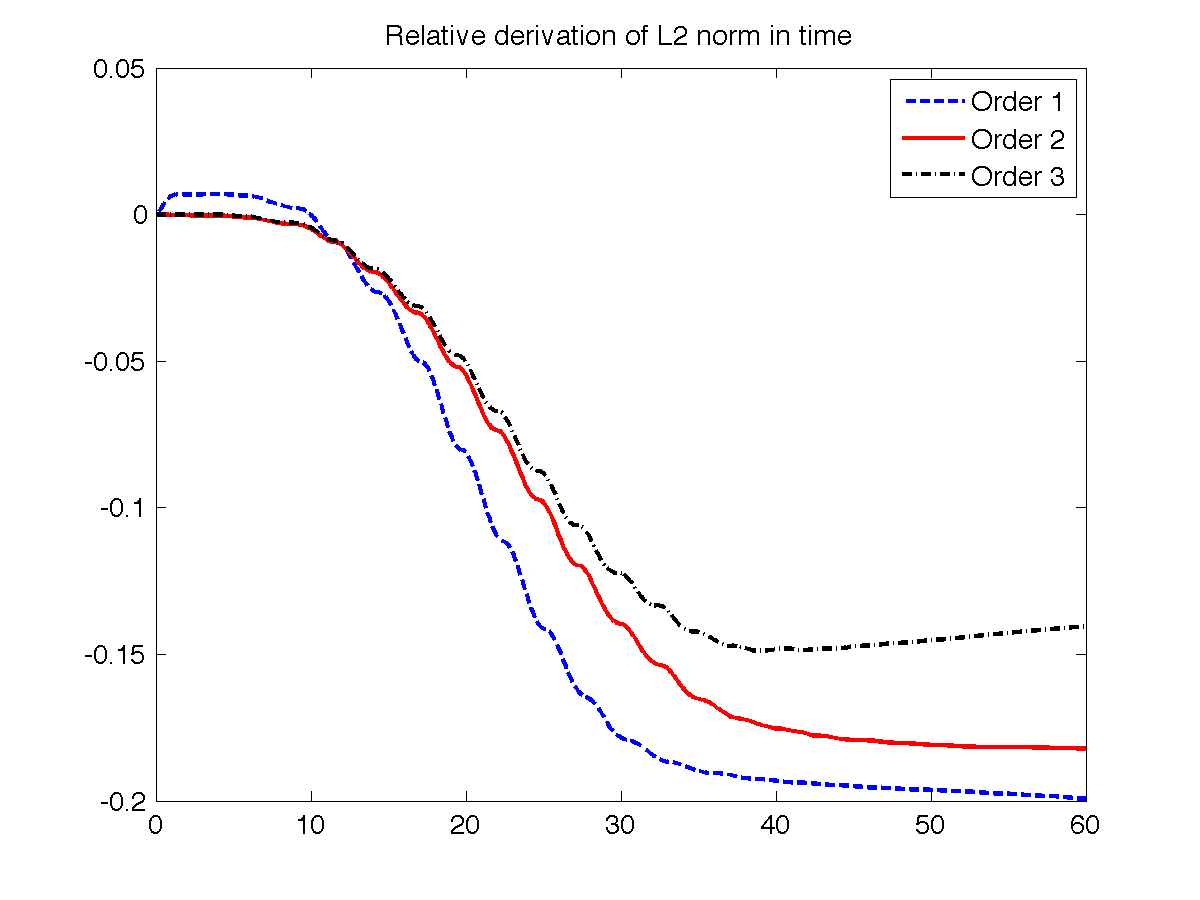}\\
\includegraphics[height=2.2in,width=3.0in]{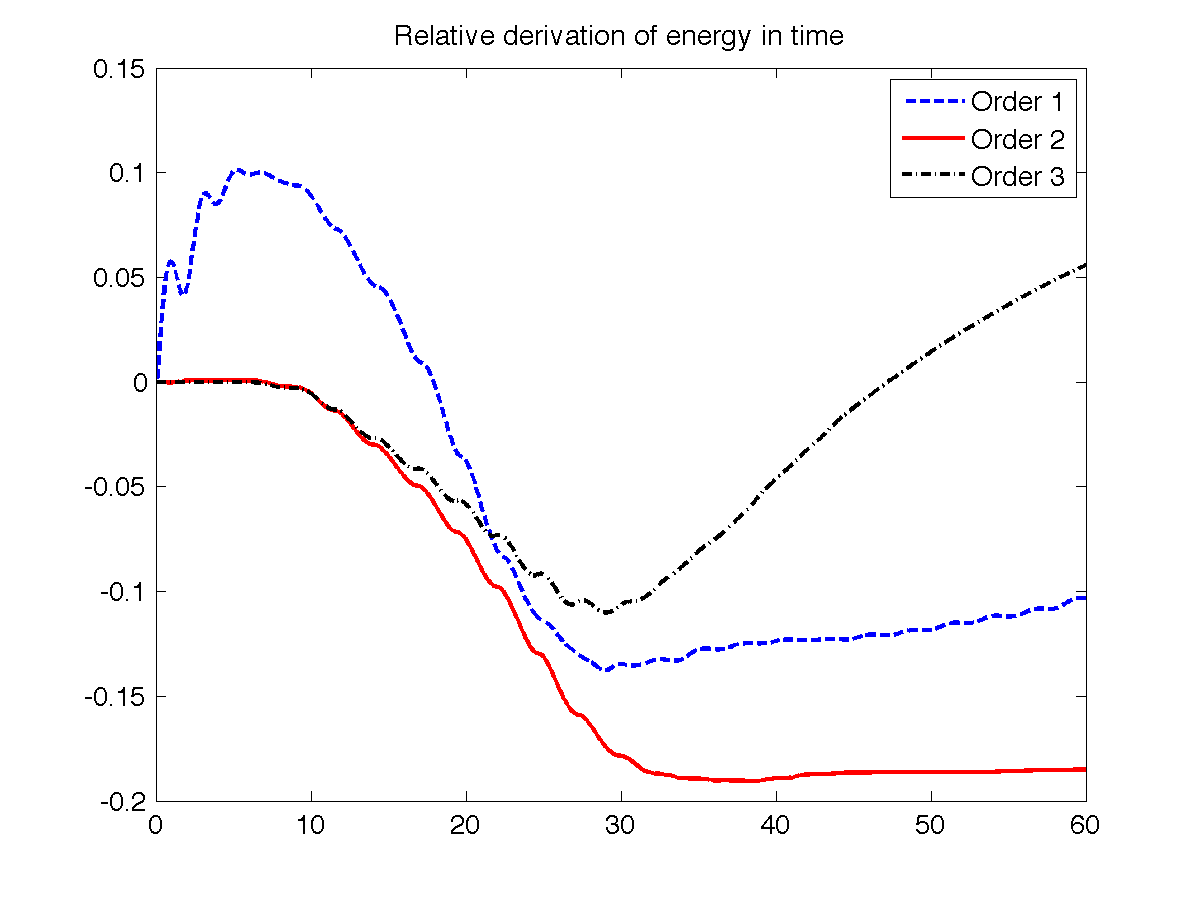}
\includegraphics[height=2.2in,width=3.0in]{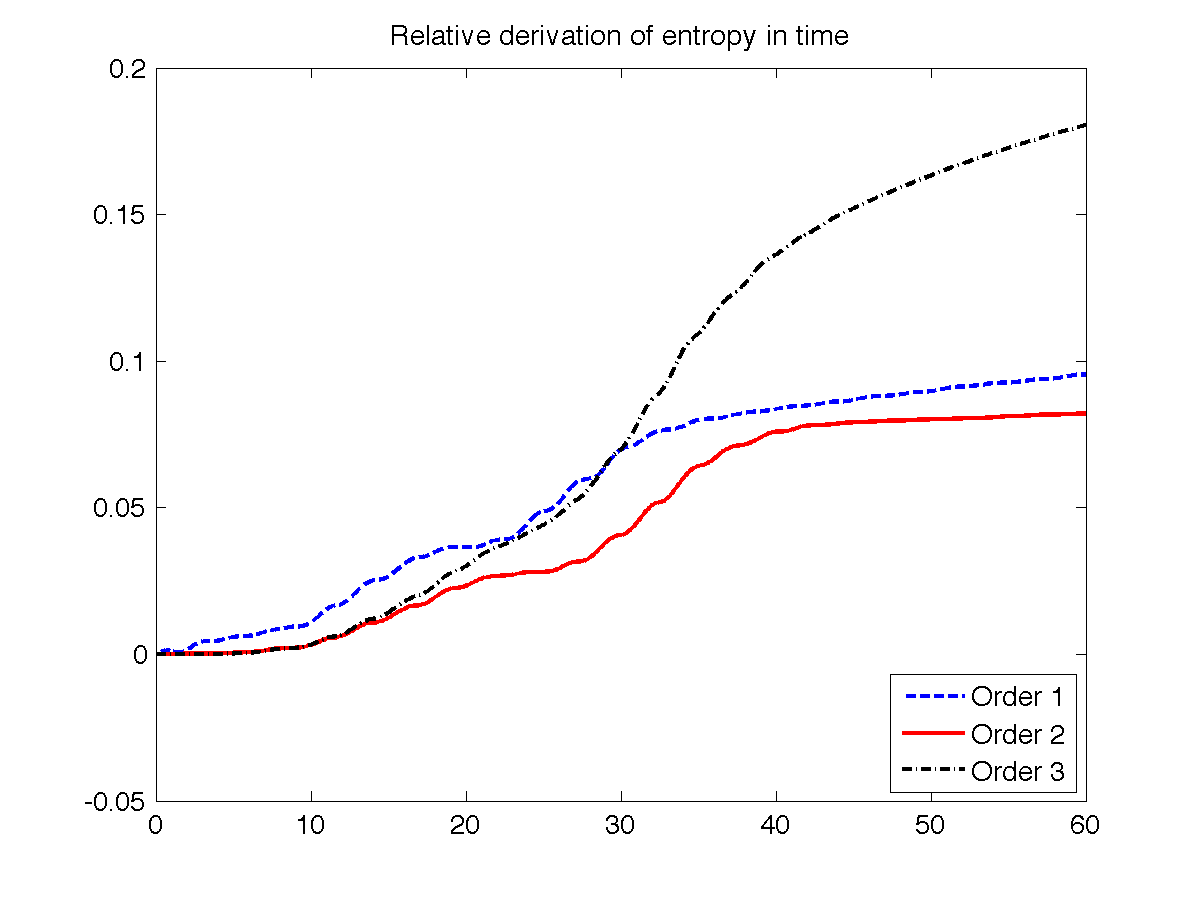}
\end{center}
\caption{Strong Landau damping. The SL WENO scheme with sixth order WENO interpolation in space and
various orders of temporal accuracy. 
Time evolution of the relative
deviations of discrete $L^1$ norms (upper left), $L^2$ norms,
kinetic energy norms (lower left) and entropy (lower right).}
\label{fig407}
\end{figure}
\end{exa}

\begin{exa} Consider the symmetric two stream
instability 
\cite{banks2010new}, 
with the initial
condition
\begin{equation}
f(x,v,t=0)=\frac{1}{\sqrt{8\pi}v_{th}}\left[\exp\left(-\frac{(v-u)^2}{2v_{th}^2}\right)+\exp\left(-\frac{(v+u)^2}{2v_{th}^2}\right)\right]\big
(1+0.0005\cos(kx)\big )
\end{equation}
with $u=5\sqrt{3}/4$, $v_{th}=0.5$ and $k=0.2$. The background ion
distribution function is fixed, uniform and chosen so that the total
net charge density for the system is zero.
Figure~\ref{fig: 2stream2_E} plots the evolution of electric fields for the proposed scheme
benchmarked with a reference rate from linear theory $\gamma = \frac{1}{\sqrt{8}}$, see \cite{banks2010new}. 
Theoretical consistent results are observed. 
Time evolution of discrete $L^1$ norm, $L^2$ norm, kinetic energy and entropy of
schemes with different temporal orders are reported in Figure
\ref{fig: 2stream2_norm}. Again, higher order schemes in general perform better
in preserving the conserved physical quantities than low order ones.
In Figures \ref{fig: 2stream2}, we report numerical
solutions from the SL WENO schemes with various temporal accuracy in approximating the distribution solution $f$. 
It can be observed that, with the same time step size, the higher order schemes (e.g. second and third order ones) perform better than a first order one. 

\begin{figure}
\begin{center}
\includegraphics[height=2.2in,width=3.0in]{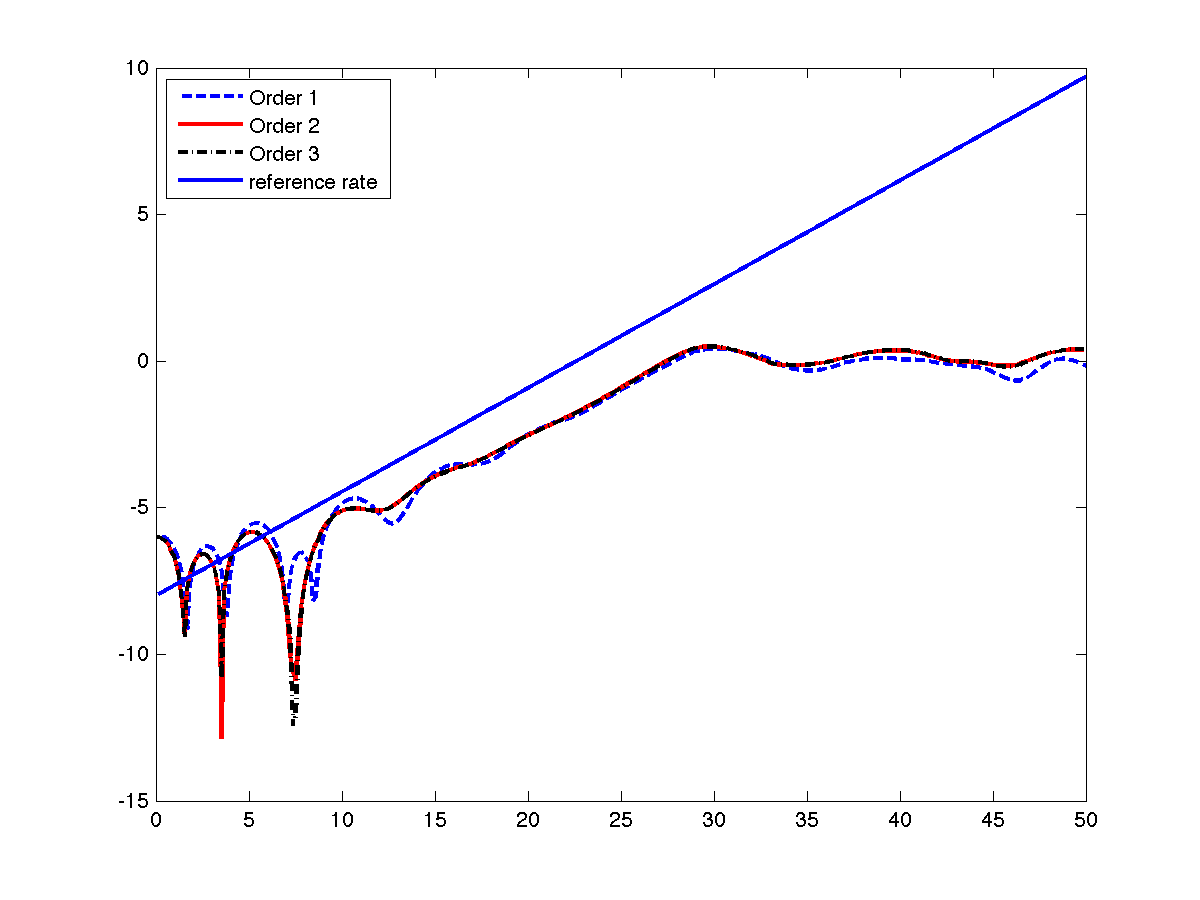}
\end{center}
\caption{Symmetric two stream instability: time evolution of electric field in
$L^2$ norm. The SL WENO scheme with sixth order WENO interpolation in space and various orders of temporal accuracy. } 
\label{fig: 2stream2_E}.
\end{figure}

\begin{figure}[htb]
\begin{center}
\includegraphics[height=2.2in,width=3.0in]{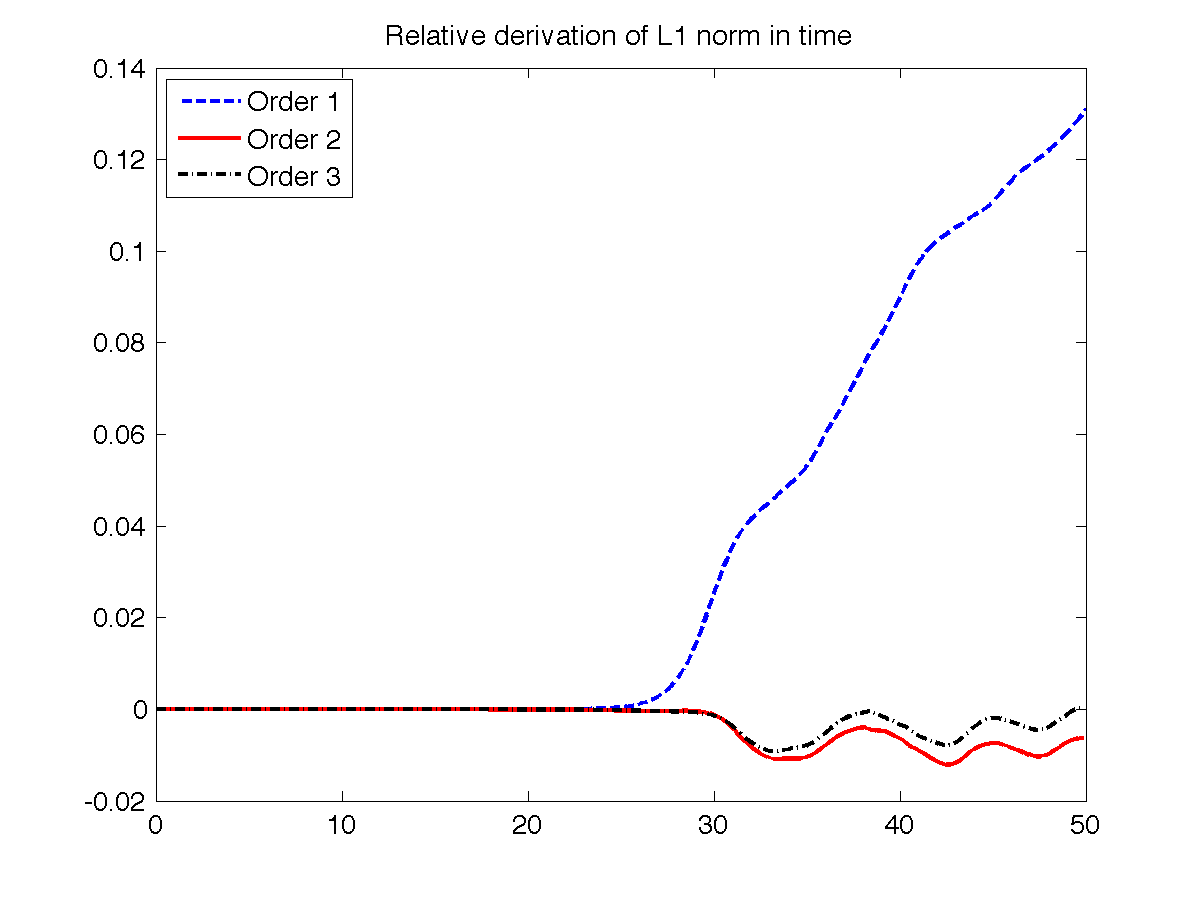}
\includegraphics[height=2.2in,width=3.0in]{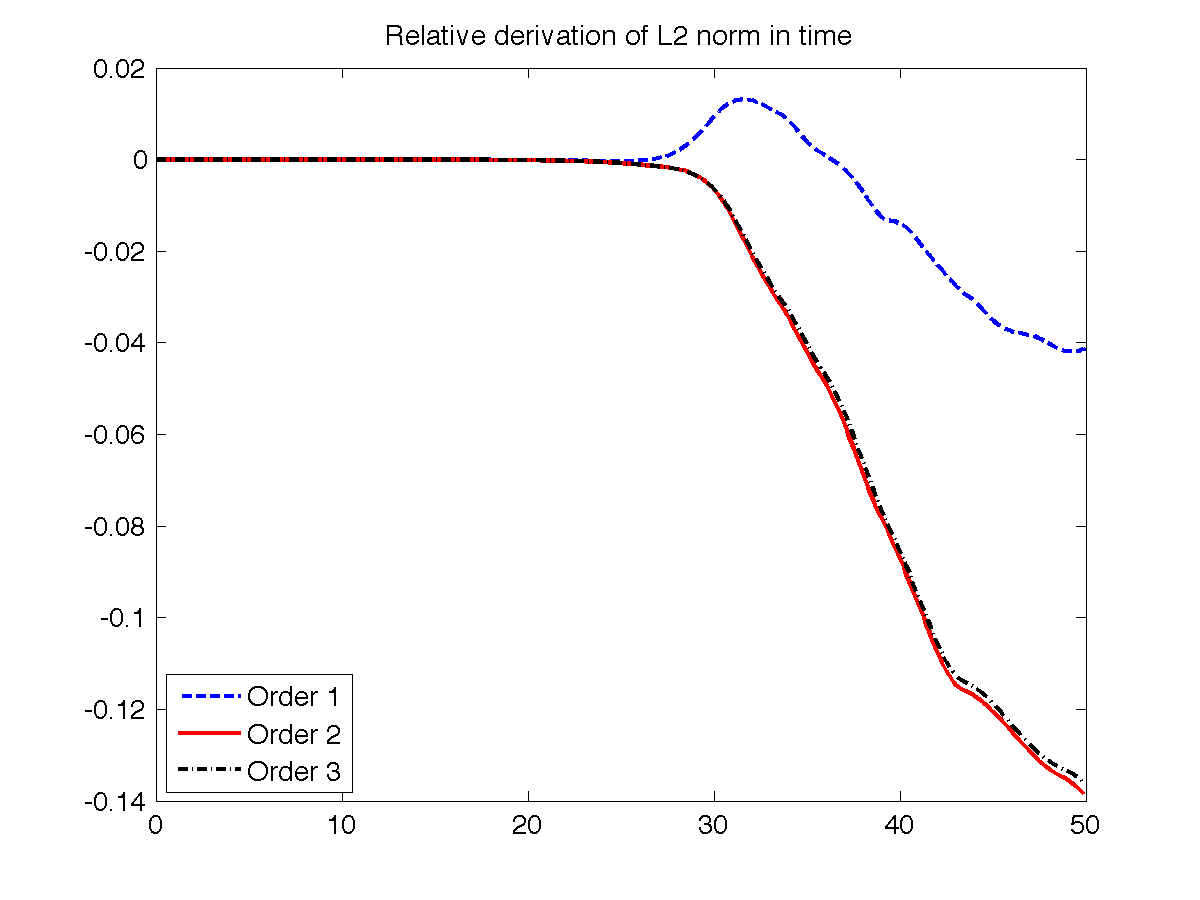}\\
\includegraphics[height=2.2in,width=3.0in]{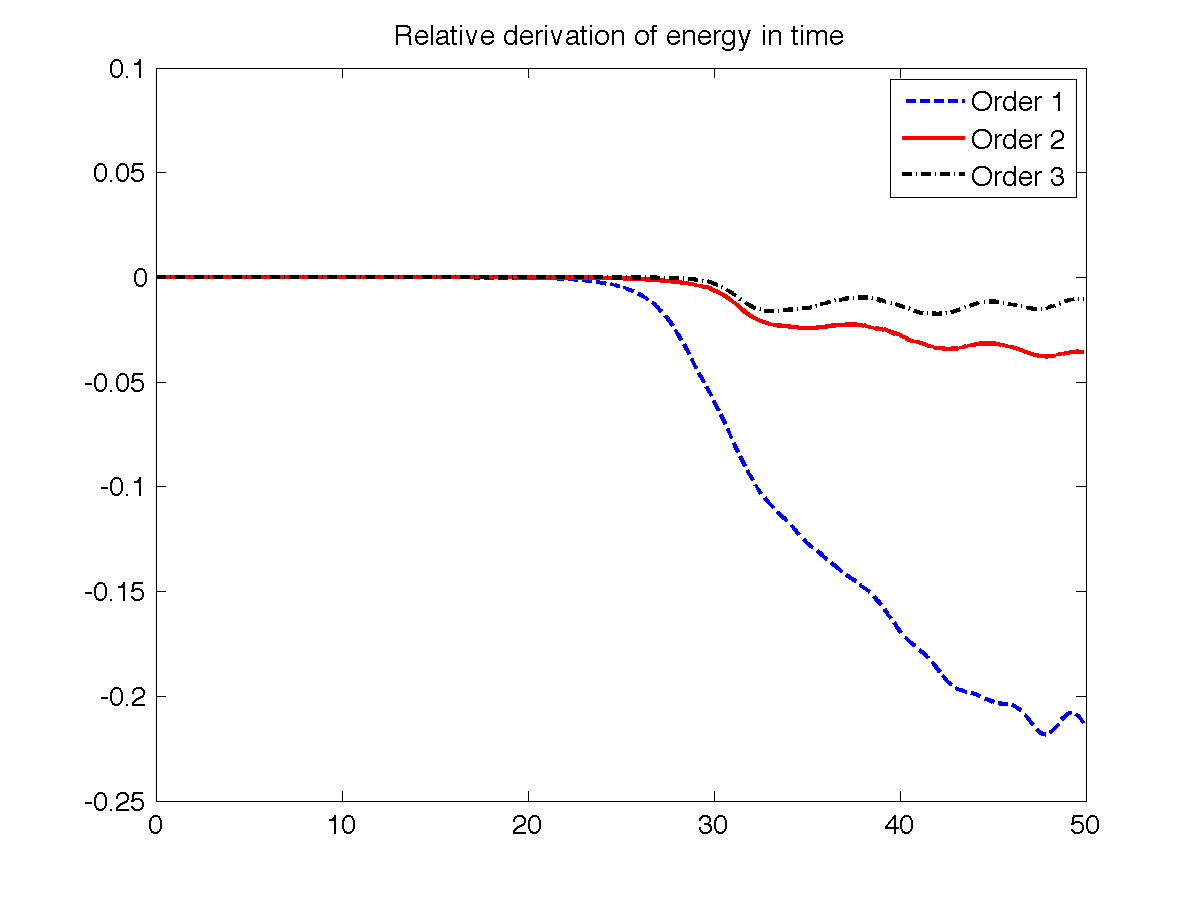}
\includegraphics[height=2.2in,width=3.0in]{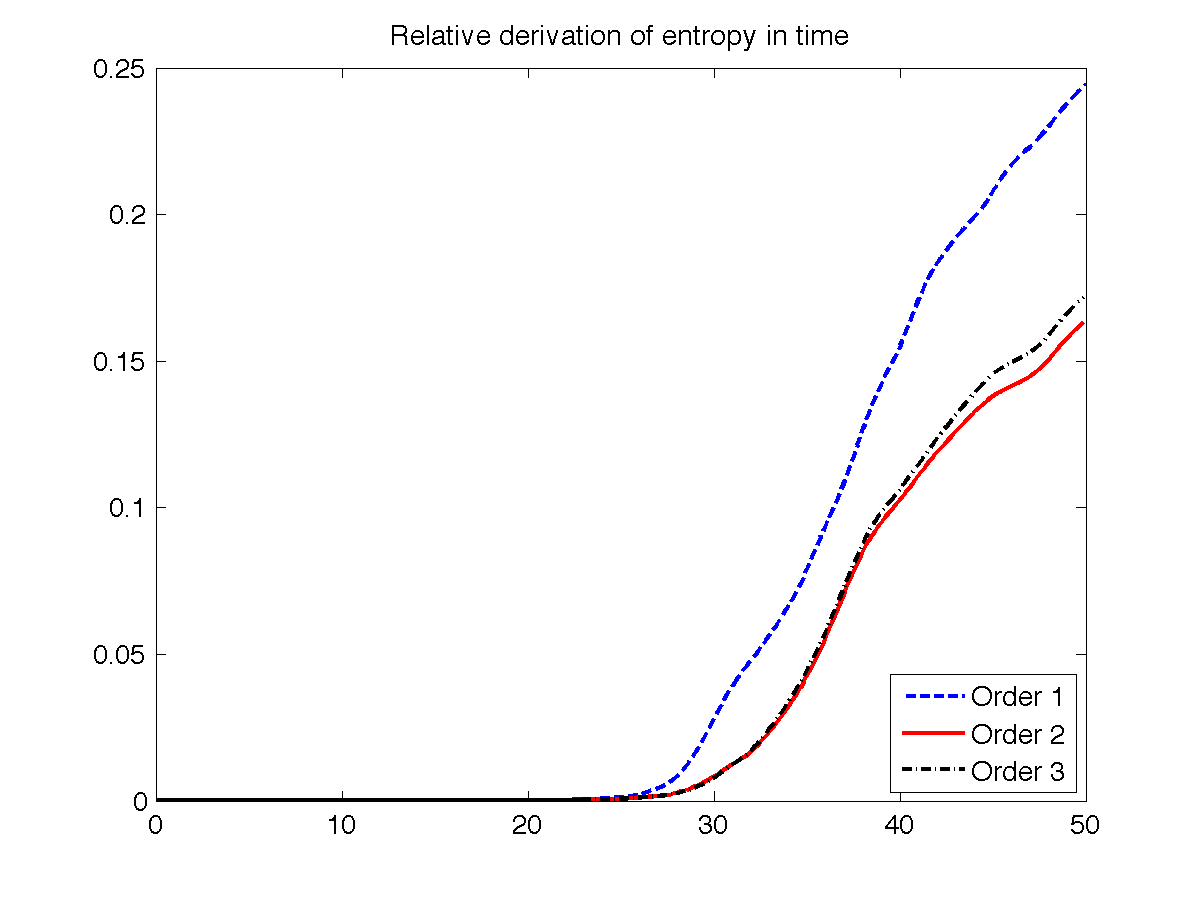}
\end{center}
\caption{Two stream instability. The SL WENO scheme with
sixth order WENO interpolation in space and various orders of temporal accuracy. 
Time evolution of the relative
deviations of discrete $L^1$ norms (upper left), $L^2$ norms,
kinetic energy norms (lower left) and entropy (lower right).}
\label{fig: 2stream2_norm}
\end{figure}

\begin{figure}[htb]
\begin{center}
\includegraphics[height=3.2in,width=3.0in]{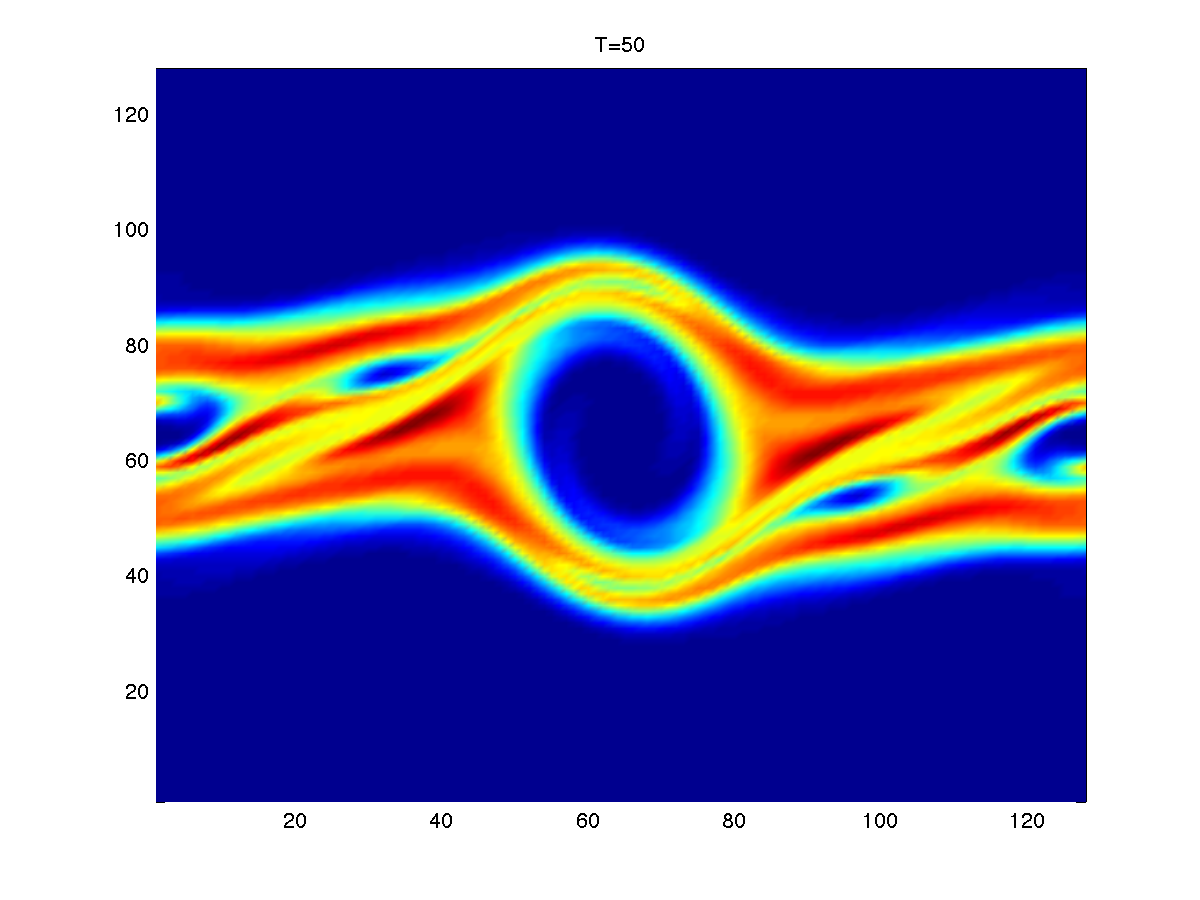}
\includegraphics[height=3.2in,width=3.0in]{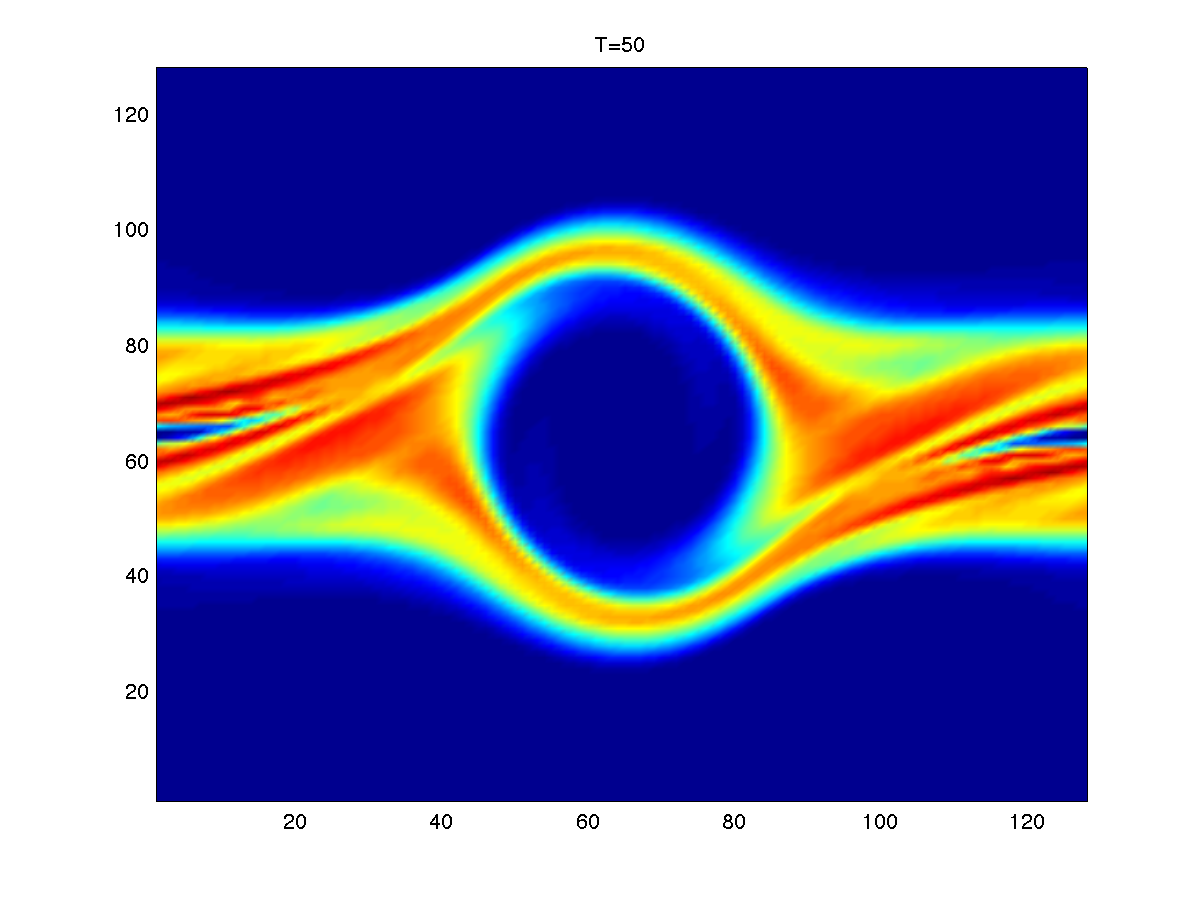}
\\
\includegraphics[height=3.2in,width=3.0in]{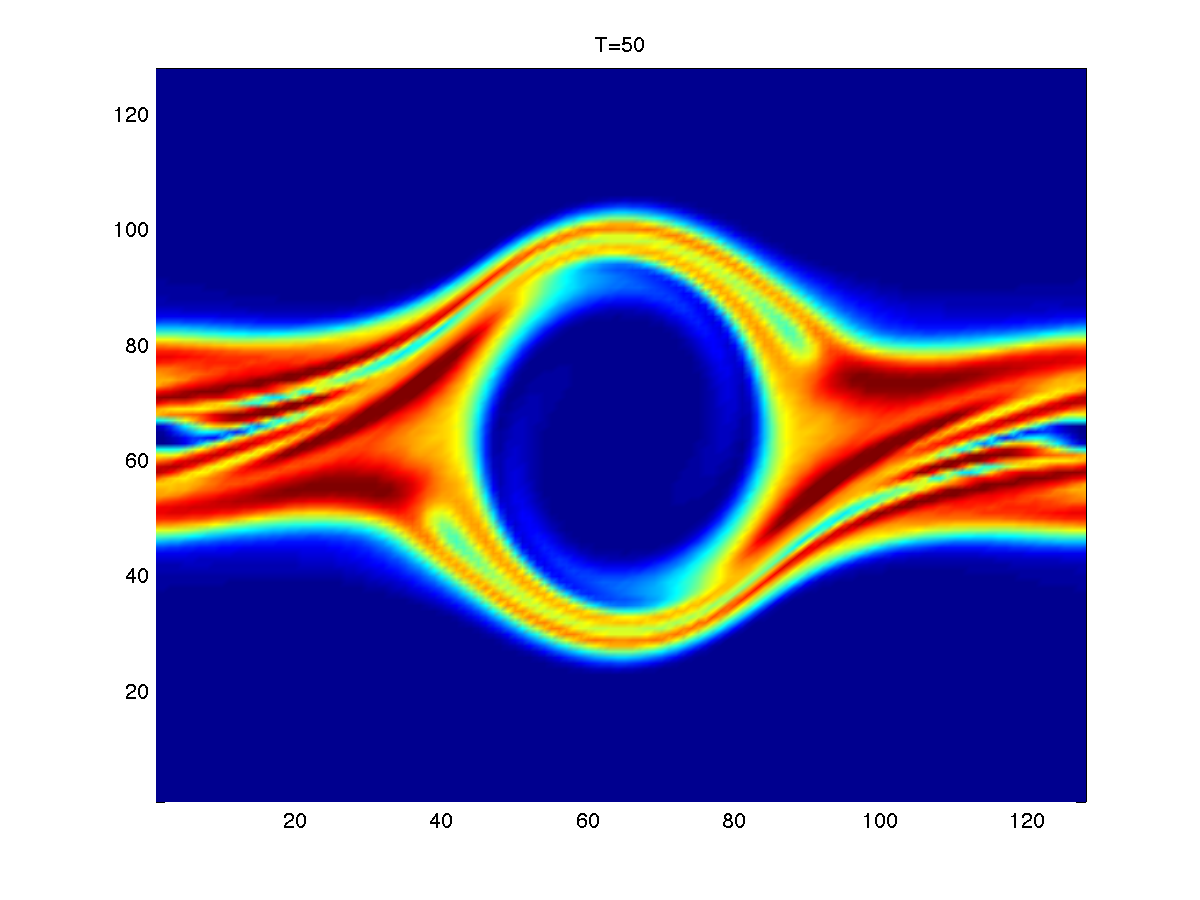}
\includegraphics[height=3.2in,width=3.0in]{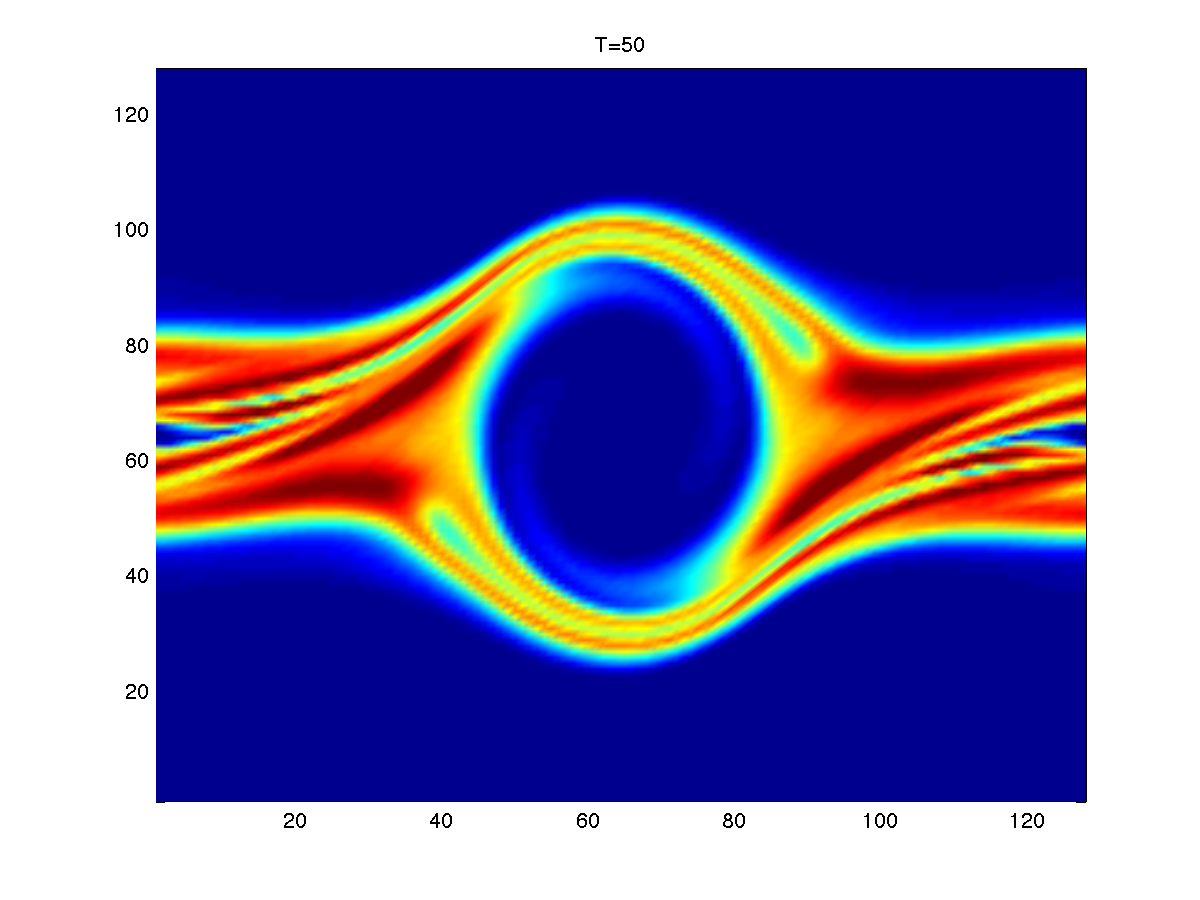}
\end{center}
\caption{Symmetric two stream instability: $T=50$. Results from schemes with first order temporal accuracy with $CFL=5$ (upper left), $CFL=0.1$ (upper right).
Results from schemes with second order temporal accuracy (lower left) and third order temporal accuracy (lower right) and $CFL=5$.}
\label{fig: 2stream2}
\end{figure}

\end{exa}

\section{Conclusion}
\label{sec5}
\setcounter{equation}{0}
\setcounter{figure}{0}
\setcounter{table}{0}

In this paper, we propose a systematical way of tracing characteristics for a one-dimensional in space and one-dimensional in velocity Vlasov-Poisson system with high order temporal accuracy. Based on such mechanism, a finite difference grid-based semi-Lagrangian approach coupled with WENO interpolation is proposed to evolve the system. It is numerically demonstrated that schemes with higher order of temporal accuracy perform better in many aspects than the first order one. Designing mass conservative semi-Lagrangian schemes, yet not subject to time step constraints, is considered to be challenging and is subject to future research investigations.

\bibliographystyle{siam}
\bibliography{refer}

\end{document}